\UseRawInputEncoding
\documentclass[11pt]{article}
\usepackage{amsmath}

\usepackage{amsmath,amssymb,amsthm,graphicx,mathrsfs}
\usepackage[numbers,sort&compress]{natbib}
\usepackage{color}
\usepackage[colorlinks,bookmarksopen,bookmarksnumbered,citecolor=black, linkcolor=blue, urlcolor=black]{hyperref}
\numberwithin{equation}{section}
\newcommand{\beq}{\begin{equation}}
\newcommand{\enq}{\end{equation}}

\newtheorem{Theorem}{Theorem}[section]
\newtheorem{Lemma}[Theorem]{Lemma}

\newtheorem{Definition}[Theorem]{Definition}
\newtheorem{Remark}[Theorem]{Remark}

\newcommand{\benu}{\begin{enumerate}}
\newcommand{\beqa}{\begin{eqnarray}}
\newcommand{\beqan}{\begin{eqnarray*}}
\newcommand{\eay}{\end{array}}
\newcommand{\edm}{\end{displaymath}}
\newcommand{\eenu}{\end{enumerate}}
\newcommand{\eeq}{\end{equation}}
\newcommand{\eeqa}{\end{eqnarray}}
\newcommand{\eeqan}{\end{eqnarray*}}

\newcommand{\br}{\begin{Remark}}
\newcommand{\er}{\end{Remark}}

\newcommand{\bqa}{\begin{eqnarray}}
\newcommand{\eqa}{\end{eqnarray}}
\newcommand{\bqw}{\begin{eqnarray*}}
\newcommand{\eqw}{\end{eqnarray*}}

\newcommand{\non}{\nonumber}
\newcommand{\bea}{\begin{array}{cc}}
\newcommand{\ena}{\end{array}}

\allowdisplaybreaks[4]
\begin{document}
\begin{center}

{\large \bf Existence and upper semicontinuity of pullback attractors for Kirchhoff wave equations in time-dependent spaces}\\

\vspace{0.20in}Bin Yang$^{1}$ $\ $ Yuming Qin $^{2,\ast}$ $\ $ Alain Miranville $^{3}$ $\ $ Ke Wang $^{2}$\\
\end{center}
$^{1}$ College of Information Science and Technology, Donghua University, Shanghai, 201620, P. R. China.\\
$^{2}$ Department of  Mathematics, Institute for Nonlinear Science, Donghua University, Shanghai, 201620, P. R. China.\\
$^{3}$ Laboratoire de Math\'ematiques et Applications, Universit\'e de Poitiers, UMR CNRS 7348-SP2MI, Boulevard Marie et Pierre Curie-T\'el\'eport 2, F-86962, Chasseneuil Futuroscope Cedex, France.\\
 \vspace{3mm}

\begin{abstract}
In this paper, we shall investigate the existence and upper semicontinuity of pullback attractors for non-autonomous Kirchhoff wave equations with a strong damping in the time-dependent space $X_t$. After deriving the existence and uniqueness of solutions by the Faedo-Galerkin approximation method, we establish the existence of pullback attractors. Later on, we prove the upper semicontinuity of pullback attractors between the Kirchhoff-type wave equations with $\delta \geq 0$ and the conventional wave equations with $\delta=0$ by a series of complex energy estimates.
\end{abstract}

\hspace{4mm}{\bf Keywords:} Kirchhoff wave equations; strong damping; pullback attractors; upper semicontinuity; the time-dependent space.

\hspace{4mm}{\bf 2020 MSC:} 35B40, 35B41, 35L05.

\section{Introduction}
\setcounter{equation}{0}

\let\thefootnote\relax\footnote{*Corresponding author: yuming@dhu.edu.cn}
\let\thefootnote\relax\footnote{\footnotesize E-mails: binyangdhu@163.com, Alain.Miranville@math.univ-poitiers.fr, kwang@dhu.edu.cn}
\quad
In this paper, we will consider the following non-autonomous Kirchhoff wave equations
\begin{equation}
\left\{\begin{array}{ll}
\varepsilon(t)\partial_{tt} u -(1+\delta \|\nabla u\|^2) \Delta u-\Delta \partial_{t}u+\lambda u=g(u)+h(x,t) & in\enspace \Omega \times (\tau +\infty), \\
u(x,t)=0 & on \enspace\partial \Omega\times(\tau +\infty), \\
u(x, \tau)=u_{\tau}^0, \,\,\partial_{t}u(x,\tau)=u_{\tau}^1, &\, x \in \Omega,
\end{array}\right.\label{1.1-6}
\end{equation}
where $\Omega \subset \mathbb R^{n}$ $(n \geq 3)$ is a bounded domain with smooth boundary $\partial\Omega$, $\varepsilon(t)$ is a time-dependent function, $-\Delta u_{t}$ is a strong damping, $g(u)$ is a nonlinear function, $h(x,t)$ is an external force and $\delta, \lambda \geq 0$ are constants.

Suppose the time-dependent function $\varepsilon(t) \in C^{1}(\mathbb{R})$ is decreasing and satisfies
\begin{equation}
\lim _{t \rightarrow+\infty} \varepsilon(t)=\alpha_{\epsilon} \geq 1,
\label{1.2-6}
\end{equation}
and there exists a constant $L\geq\alpha_{\epsilon}$ such that
\begin{equation}
\sup _{t \in \mathbb{R}}(|\varepsilon(t)|+|\varepsilon^{\prime}(t)|) \leq L.
\label{1.3-6}
\end{equation}

Moreover, assume $g \in {C}^1(\mathbb R;\mathbb R)$ with $g(0)=0$ and satisfies
\begin{equation}\label{1.4-6}
\left|g^{\prime}(u)\right| \leq C\left(1+|u|^{p}\right),
\end{equation}
\begin{equation}\label{1.5-6}
g^{\prime}(u) \leq k,
\end{equation}
\begin{equation}\label{1.6-6}
\limsup _{|u| \rightarrow \infty} \frac{u g(u)-\gamma G(u)}{u^2} \leq 0
\end{equation}
and
\begin{equation}\label{1.7-6}
\limsup _{|u| \rightarrow \infty} \frac{G(u)}{u^2} \leq 0,
\end{equation}
where $p=\frac{4}{n-2}$, $C,k,\gamma>0$ are constants and $G(u)=\int_0^u g(s) d s$.

Besides, let $h(x, t), \partial_t h(x, t) \in L_{\text {loc }}^2\left(\mathbb{R} ; L^2(\Omega)\right)$ and there exists a parameter $\sigma>0$ such that
\begin{equation}
\int_{-\infty}^t e^{\sigma s}\|h(x, s)\|^2 d s<+\infty, \quad \forall \, t \in \mathbb{R}.
\label{1.8-6}
\end{equation}

Our phase space is the time-dependent space $X_t=H_{0}^1(\Omega)\times L^{2}(\Omega)$, which is equipped with the following norm
\begin{equation}
\|(u_0, u_1)\|_{X_t}=\|\nabla u_0\|^2+\varepsilon(t)\|u_1\|^2.
\label{1.9-6}
\end{equation}

Additionally, the norm of $X_t$ can be equivalently written as
\begin{equation}
\|(u, \partial_t u)\|^2=\|\nabla u\|^2+\varepsilon(t)\|\partial_t u\|^2.
\label{1.11-6}
\end{equation}

In 1883, Kirchhoff \cite{k.5} first proposed problem $u_{tt}-\left(1+\epsilon\|\nabla u\|^2\right) \Delta u+f(u)=g(x)$ for any $\epsilon >0$.
Subsequently, scholars discovered that the Kirchhoff-type wave equations play a crucial role in seismic imaging, oil and gas exploration, and other fields. Therefore, many mathematicians have devoted themselves to studying attractors for such equations (see \cite{c.5, l.5, ljs.5, mz.5}).
It is noted that the Kirchhoff term is non-local, which brings some difficulties to derive the well-posedness of solutions. In addition, if a Kirchhoff wave equation is non-autonomous and contains a damping, the difficulties will be greatly increased. So far, few results have been obtained about the attractors for the non-autonomous Kirchhoff wave equations with a damping.

Fan and Zhou \cite{fz.6} proved the existence of compact kernel sections for the process generated by problem $u_{t t}-\alpha \Delta u_t-\left(\beta+\gamma\left(\int_{\Omega}|\nabla u|^2 \mathrm{~d} x\right)^\rho\right) \Delta u+h\left(u_t\right)+f(u, t)=g(x, t)$ in $H_0^1(\Omega) \times L^2(\Omega)$, where $\alpha, \beta>0$, $\rho>-1$ and $\gamma \geq 0$ are constants.
Moreover, \cite{fz.6} derived a precise estimate of upper bound of Hausdorff dimension of kernel sections, which decreases as the strong damping grows for the large strong damping, particularly in the autonomous case. Then Wang and Zhong \cite{wz.5} considered the upper semicontinuity of pullback attractors for problem $u_{tt}-\Delta u_t-\left(1+\epsilon_1 \|\nabla u\|^2\right) \Delta u+$ $f(u)=g(x, t)$ with $\epsilon_1 >0$ in $H_0^1(\Omega) \times L^2(\Omega)$.
Furthermore, Yang and Li \cite{YL.6} concerned with the existence and stability of pullback exponential attractors for problem $u_{t t}-M\left(\|\nabla u\|^2\right) \Delta u+(-\Delta)^{\alpha_1} u_t+f(u)=g(x, t)$ in $H^{1+\alpha_1}(\Omega)\times H^{\alpha_1}(\Omega)$ for $\alpha_1 \in(1 / 2,1)$.
Later on, Li and Yang \cite{ly.5} proved the robustness of pullback attractors and pullback exponential
attractors for problem $u_{t t}-\left(1+\epsilon_2\|\nabla u\|^2\right) \Delta u-\sigma\left(\|\nabla u\|^2\right) \Delta u_t+f(u)=g(x, t)$ in $\left(H_0^1(\Omega) \cap L^{p+1}(\Omega)\right) \times L^2(\Omega)$, where $\epsilon_2 \in[0,1]$ is an extensibility parameter.

Yang and Li \cite{LYUPP.6} obtained the process generated by problem $u_{t t}-M\left(\|\nabla u\|^2\right) \Delta u+(-\Delta)^{\alpha_2} u_t+f(u)=g(x, t)$ has pullback attractors for $\alpha_2 \in(1 / 2,1)$ and the family of pullback attractors is upper semicontinuous in $H_0^1 (\Omega)\cap L^{p+1} (\Omega)\times L^2(\Omega)$ when the nonlinearity $f(u)$ is of supercritical growth $p: 1 \leq p<p_{\alpha_2} \equiv \frac{N+4 \alpha_2}{(N-4 \alpha_2)^{+}}$ with $N\geq3$.
In addition, Ma, Wang and Xie \cite{MWX2021-Y0} verified the existence of pullback attractors for problem $u_{t t}-\Delta u_t-\phi\left(\|\nabla u\|^2\right) \Delta u+f(u)=h(x, t)$ in $H_0^1(\Omega) \times L^2(\Omega)$.
Furthermore, Li, Yang and Feng \cite{lyf.5} established the existence and continuity of uniform attractors for problem $u_{t t}-\left(1+\epsilon_3\|\nabla u\|^2\right) \Delta u-\Delta u_t+f(u)=g(x, t)$ in $H_0^1(\Omega) \cap L^{p+1}(\Omega) \times L^2(\Omega)$ when $f(u)$ satisfies optimal supercritical growth condition $p: \frac{N+2}{N-2}=p^*<p<p^{* *}=\frac{N+4}{(N-4)^{+}}$ with $N\geq3$ and $\epsilon_3 \in[0,1]$.

Moreover, we will state the main difficulties and innovations of this paper.

(i) It is well-known that Kirchhoff wave equations have a wide range of applications in actual life. However, since each non-autonomous Kirchhoff wave equation contains a non-local Kirchhoff term and a time-dependent external force function, which bring lots of computational difficulties to investigate the long-time behavior of solutions. Therefore, there are only a dozen papers about the existence and upper semicontinuity of pullback attractors for Kirchhoff wave equations. Meanwhile, no paper has studied pullback attractors of Kirchhoff wave equations in time-dependent space $X_t=H_{0}^1(\Omega)\times L^{2}(\Omega)$, thus our paper is a new attempt.

(ii) In addition, the strong damping $-\Delta \partial_{t}u$ makes the application of problem \eqref{1.1-6} more comprehensive.
But $-\Delta \partial_{t}u$ directly makes it necessary to overcome more difficulties in calculations and estimates. We deal with these difficulties by undertaking a large number of meticulous energy estimates and by adjusting the range of parameters in test functions.

(iii) When $\delta=0$, problem \eqref{1.1-6} is an ordinary wave equation, in which case we obtain that there exists a pullback attractor $A_0$. Moreover, when $\delta\geq0$, problem \eqref{1.1-6} is a Kirchhoff wave equation with  a strong damping, in this case we derive that there exists a pullback attractor $A_\delta$. We verify the upper semicontinuity of $A_0$ and $A_\delta$ in Section 4, which is helpful for learning the relationship between Kirchhoff wave equations and ordinary wave equations.

This paper is organized as follows. In Section 2, we recall some definitions and lemmas about pullback attractors. Furthermore, by a series of elaborate energy estimates, in Section 3 we derive the existence of  pullback attractors for the process generated by problem \eqref{1.1-6} in the time-dependent space $X_t$. In Section 4, we establish the upper semicontinuity of pullback attractors between the non-autonomous Kirchhoff wave equations \eqref{1.1-6} with $\delta \geq 0$ and the conventional wave equations with $\delta=0$ in $X_t$.

\section{Preliminaries}
$\ \quad$ In this section, we will recall some necessary definitions and lemmas about pullback attractors and upper semicontinuity of pullback attractors in detail.

Suppose $X$ is a Banach space or a closed subset of a Banach space, which is endowed with the norm $\|\cdot\|_{X}$.
\begin{Definition} {\rm(\cite{ConPT.6, r})}
A process or a two-parameter semigroup on $X$ is the family $\{U(t, \tau)\}_{t \geq \tau}$ of mapping $U(t, \tau): X \rightarrow X$ satisfying that $U(\tau, \tau)=Id$ is the identity operator and $U(t, s) U(s, \tau)=U(t, \tau)$ for any $t \geq s \geq \tau$.
\label{def2.1-6}
\end{Definition}

\begin{Definition} {\rm(\cite{ConPT.6, r})}
The Hausdorff semi-distance of two nonempty  sets $A_1, A_2 \subset X$ is defined by
$$
dist_{X}(A_1, A_2)=\sup _{x_1 \in A_1} \inf _{x_2 \in A_2}\|x_1-x_2\|_{X} \, .
$$
\end{Definition}

\begin{Definition} {\rm(\cite{bv1992,psz,zs2})}
A family ${\mathcal D=\{ D(t)\}}_{t \in \mathbb{R}}$ is pullback absorbing for the process $U(\cdot,\cdot)$ in $X$, if for any $t \in \mathbb{R}$ and any bounded subset $D \subset X$, there exists a constant $T=T(t, D)>0$ such that $U(t, t-\tau) D \subset D(t)$ for any $\tau \geq T.$
\label{def2.2-6}
\end{Definition}

\begin{Definition} {\rm(\cite{bv1992,psz,zs2})}
A process $U(\cdot, \cdot)$ is pullback $\mathcal{D}$-asymptotically compact in $X$, if for any $t \in \mathbb{R}$, any sequence $\tau_i \rightarrow -\infty$ and $x_i \in D\left(t-\tau_i\right)$, the sequence $\left\{U\left(t, t-\tau_i\right) x_i\right\}_{i \in \mathbb{N}}$ is relatively compact in $X$.
\label{def2.3-6}
\end{Definition}

\begin{Definition} {\rm(\cite{bv1992,psz,zs2})}
A family of compact sets $\mathcal{A}=\{A(t)\}_{t \in \mathbb{R}}$ is called a pullback attractor for the process $U(\cdot, \cdot)$, if it fulfills the following properties

\noindent(i) $\mathcal{A}$ is invariant, i.e., $U(t, \tau) A(\tau)=A(t)$ for any $\tau \leq t$;

\noindent(ii) $\mathcal{A}$ is pullback attracting, i.e.,
$$\lim\limits_{\tau \rightarrow -\infty} \operatorname{dist}_X(U(t, t-\tau) D, A(t))=0$$
for any bounded $D \subset X$.
\label{def2.4-6}
\end{Definition}

\begin{Lemma}{\rm(\cite{clr.6})}\label{lem2.5-6}
Suppose the family $\mathcal{D}=\{D(t)\}_{t \in \mathbb{R}}$ is pullback absorbing and the process $U(\cdot, \cdot)$ is pullback $\mathcal{D}$-asymptotically compact in $X$, if the family $\mathcal{A}=\{A(t)\}_{t \in \mathbb{R}}$ satisfies
\begin{align}\label{1a-6}
A(t)=\Lambda(\mathcal{D}, t)= \mathop  \bigcap \limits_{s \ge 0} {\overline {\mathop  \bigcup \limits_{\tau  \ge s} U(t,t - \tau )D(t - \tau )} ^X},\quad \forall \, t \in \mathbb R,
\end{align}
then $\mathcal A$ is a pullback attractor for $U(\cdot, \cdot)$ in $X$. In addition, if for any $t \in \mathbb{R}$, there exists a $T=T(t)>0$ such that
\begin{equation}\label{1b-6}
U(t, t-\tau) D(t-\tau) \subset D(t), \quad \forall \, \tau>T,
\end{equation}
then
\begin{equation}\label{1c-6}
\lim _{\tau \rightarrow -\infty} \operatorname{dist}_X(U(t, t-\tau) D(t-\tau), A(t))=0.
\end{equation}
\end{Lemma}

\begin{Definition} {\rm(\cite{bv1992,psz,zs2})}\label{def2.6-6}
Suppose the set $D \subset X$ is bounded, then $\Psi(\cdot, \cdot): D \times D$ is called a  contractive function, if for any sequence $\left\{x_i\right\}_{i \in \mathbb{N}} \subset D$, there exists a subsequence $\left\{x_{i_k}\right\}_{k \in \mathbb{N}} \subset$ $\left\{x_i\right\}_{i \in \mathbb{N}}$ such that
$$
\lim _{k \rightarrow +\infty} \lim _{l \rightarrow +\infty} \Psi\left(x_{i_k}, x_{i_l}\right)=0.
$$
\end{Definition}

\begin{Lemma}{\rm(\cite{wyh1.6})}\label{lem2.7-6}
Suppose the family $\mathcal{D}=\{D(t)\}_{t \in \mathbb{R}} \subset X$ and satisfies \eqref{1b-6}, and for any constant $\zeta>0$ and any $t \in \mathbb{R}$, there exists a time $\tau=\tau(t, \zeta,\mathcal D)>0$ and a contractive function $\Psi_{t-\tau}(\cdot, \cdot): D(t-\tau) \times D(t-\tau)$ depending on $t-\tau$ such that
$$
\|U(t, t-\tau) x-U(t, t-\tau) y\|_X \leq \zeta+\Psi_{t-\tau}(x, y), \quad \forall \, x, y \in D(t-\tau),
$$
then $U(\cdot, \cdot)$ is pullback $\mathcal{D}$-asymptotically compact in $X$.
\end{Lemma}

Below we will introduce the correlation lemma, which is used to prove the process $U(\cdot, \cdot)$ is pullback $\mathcal{D}$-asymptotically compact in $X$, as well as the lemmas to derive the upper semicontinuity of pullback attractors, which have been systematically proved in \cite{wz.5}.

Let $X_a$ and $X_b$ be two Banach spaces, which are endowed with the norms $\| \cdot \|_{X_a}$ and $\| \cdot \|_{X_b}$, respectively. Then assume $X=X_a \times X_b$ is equipped with the norm $\| \cdot \|_X$. In addition, suppose the process $U(\cdot,\cdot)$ defined in $X$ admits that $U(\cdot,\cdot):\left(u_0 ; u_1\right) \rightarrow\left(U_a(\cdot, \cdot) u_0 ; U_b(\cdot, \cdot) u_1\right)$, where $U_a(\cdot, \cdot)$ and $U_b(\cdot, \cdot)$ are continuous and defined in $X_a$ and $X_b$, respectively.
\begin{Lemma}{\rm(\cite{wz.5})}
Suppose the process $U_a(\cdot, \cdot)=U_{a_1}(\cdot, \cdot)+U_{a_2}(\cdot, \cdot)$ and the family $\mathcal{D}=\{D(t)\}_{t \in \mathbb{R}}\subset X$, if for any $t \in \mathbb R$ the following properties hold

\noindent(i) there exists a function $\phi: \mathbb{R} \times \mathbb{R} \rightarrow \mathbb{R}^{+}$ satisfying
$$\lim\limits _{\tau \rightarrow -\infty} \phi(t, \tau)=0$$
and
$$
\left\|\bigcup_{\left(u_0 ; u_1\right) \in D(t-\tau)} U_{a_1}(t, t-\tau) u_0\right\|_{X_a} \leq \phi(t, \tau)
$$
for any $\tau>0$;

\noindent(ii) $\bigcup\limits_{\left(u_0 ; u_1\right) \in D(t-\tau)} U_{a_2}(t, t-\tau) u_0$ is relatively compact in $X_a$ for any $\tau>0$;

\vspace{0.3cm}
\noindent(iii) $\bigcup\limits_{\left(u_0 ; u_1\right) \in D(t-\tau)} U_b(t, t-\tau) u_1$ is relatively compact in $X_b$ for any $\tau>0$;

\vspace{0.3cm}
\noindent then $U(\cdot, \cdot)$ is pullback $\mathcal{D}$-asymptotically compact in $X$.
\label{lem2.8-6}
\end{Lemma}

\begin{Lemma}{\rm(\cite{wz.5})}
Assume the process $U_{\delta}(\cdot, \cdot)$ has a family of pullback absorbing sets $\mathcal{D}_{\delta}=\left\{D_{\delta}(t)\right\}_{t \in \mathbb{R}}$ satisfying \eqref{1b-6} for any $\delta \in I$, $\mathcal{A}_{\delta}=$ $\left\{A_{\delta}(t)\right\}_{t \in \mathbb{R}}$ satisfies Lemma $\ref{lem2.5-6}$, and the following properties hold

\noindent(i) there exists a $\delta_0 \in I$ satisfying that
\begin{equation}\label{1d-6}
{\overline {\bigcup\limits_{\delta  \in I} {{D_\delta }(t)} } ^X} \subset {D_{{\delta _0}}}(t);
\end{equation}

\noindent(ii) for any sequences $\delta_i \in I$ with $\delta_i \rightarrow \delta_0 \in I$, $t_i>\alpha$ with $t_i \rightarrow t_0$, $x_i \in X$ with $x_i \rightarrow x_0$ and $\tau>0$, it follows
\begin{equation}\label{1e-6}
\lim _{i \rightarrow +\infty}\left\|U_{\delta_i}\left(t_i, \alpha-\tau\right) x_i-U_{\delta_0}\left(t_0, \alpha-\tau\right) x_0\right\|_X=0;
\end{equation}
(iii) for any sequence $\delta_i \in I$, any bounded $\left\{t_i\right\}_{i \in \mathbb{N}} \subset \mathbb{R}$, $\tau_i>0$ with $\tau_i \rightarrow -\infty$ and $x_i \in D\left(\alpha-\tau_i\right)$,
\begin{equation}\label{1f-6}
\left\{U_{\delta_i}\left(t_i, \alpha-\tau_i\right) x_i\right\}_{i \in \mathbb{N}} \text { is relatively compact in } X;
\end{equation}
then
\begin{equation}\label{1g-6}
\lim _{\delta \rightarrow \delta_0} \sup _{t \in[\alpha, \beta]} \operatorname{dist}_X\left(A_{\delta}(t), A_{\delta_0}(t)\right)=0
\end{equation}
for any $[\alpha, \beta] \subset \mathbb{R}$.
\label{lem2.9-6}
\end{Lemma}

The following lemma is used to prove Lemma \ref{lem2.9-6} (iii).
\begin{Lemma}{\rm(\cite{wz.5})}
Assume $[\alpha, \beta] \times I \subset \mathbb{R} \times \mathbb{R}$ and there exists a constant $W=W(t, T)>0$ such that the pullback absorbing family $\mathcal D_{\delta}=\left\{D_{\delta}(t)\right\}_{t \in \mathbb{R}}$ satisfying
$$U_{\delta}(t, t-\tau)D_{\delta}(t-\tau)\subset D_{\delta}(t)$$
for any $t \in \mathbb R$ and $\tau>W$. Additionally, let $U_{\delta}(\cdot, \cdot)=\left(U_a(\cdot, \cdot) ; U_{b}(\cdot, \cdot)\right)$ and $U_{a}(\cdot, \cdot)=$ $U_{a_1}(\cdot, \cdot)+U_{a_2}(\cdot, \cdot)$, if the following properties hold

\noindent(i) there exists a function $\Phi_{[\alpha, \beta]}(\cdot): \mathbb{R} \rightarrow \mathbb{R}^{+}$ satisfying
$$\lim_{\tau \rightarrow -\infty} \Phi_{[\alpha, \beta]}(\tau)=0$$
and
$$\left\|\bigcup_{t \in[\alpha, \beta]} \bigcup_{\delta \in I} \bigcup_{\left(u_0 ; u_1\right) \in D_{\delta}(\alpha-\tau)} U_{a_1}(t, \alpha-\tau) u_0\right\|_{X_a} \leq \Phi_{[\alpha, \beta]}(\tau) ;$$
\vspace{0.3cm}
\noindent (ii) $\bigcup\limits_{t \in[\alpha, \beta]} \bigcup\limits_{\delta \in I} \bigcup\limits_{\left(u_0 ; u_1\right) \in D_{\delta}(\alpha-\tau)} U_{a_2}(t, \alpha-\tau) u_0$ is relatively compact in $X_a$ for any $\tau>0$;\\
(iii) $\bigcup\limits_{t \in[\alpha, \beta]} \bigcup\limits_{\delta \in I} \bigcup\limits_{\left(u_0 ; u_1\right) \in D_{\delta}(\alpha-\tau)} U_{b}(t, \alpha-\tau) u_1$ is relatively compact in $X_b$ for any $\tau>0$;

\vspace{0.3cm}
\noindent then Lemma $\ref{lem2.9-6}$ (iii) holds.
\label{lem2.10-6}
\end{Lemma}

\section{Existence of pullback attractors $\mathcal A_\delta=\{A_\delta(t)\}_{t\in \mathbb R}$}
$\ \quad$ In this section, we first derive the existence and uniqueness of weak solutions for problem \eqref{1.1-6} in the time-dependent space $X_t$. Then through a series of elaborate and complex calculations and estimates, we establish the existence of the pullback attractors $\mathcal A_\delta=\{A_\delta(t)\}_{t\in \mathbb R}$ for the process $U_{\delta}(\cdot,\cdot)$ generated by problem \eqref{1.1-6} in $X_t$.

Using the Faedo-Galerkin approximation method, the following theorem can be easily verified.
\begin{Theorem}
Under the assumptions of $\varepsilon(t)$, $g(u)$, $h(x,t)$, $\delta$ and $\lambda$ in Section 1, let $\left(u_\tau^0, u_\tau^1\right)$ be given, then for any $\delta \geq 0$ and any $T>\tau$, there exists a unique weak solution $u$ to problem \eqref{1.1-6} satisfying
$$
(u, \partial_tu) \in C([\tau, T] ; X_t)
$$
in the time-dependent space $X_t$, which continuously depends on its initial data in $X_t$.
\label{th3.1-6}
\end{Theorem}

For any weak solution $u$ to problem \eqref{1.1-6}, let $\varphi_t=(u(t), \partial_t u(t))$, then we can define the continuous process $U_\delta(t, \tau)=X_\tau \rightarrow X_t$ generated by problem \eqref{1.1-6} as follows
$$
U_\delta(t, \tau) \varphi_\tau=\left(u(t), \partial_t u(t)\right)
$$
for any $t \geq \tau$, $\delta \geq 0$ and $\varphi_\tau \in X_\tau$.

\begin{Lemma} \label{lem3.2-6}
Under the assumptions of $\varepsilon(t)$, $g(u)$, $h(x,t)$, $\delta$ and $\lambda$ in Section 1, then it follows that
\begin{equation}\label{3.1-6}
\begin{aligned}
\|\nabla u(t)\|^2+\varepsilon(t)\left\|\partial_t u(t)\right\|^2 & \leq C e^{-\sigma_1(t-\tau)}\left(\left\|\nabla u_\tau^0\right\|^2+\varepsilon(\tau)\left\|u_\tau^1\right\|^2+\left\|\nabla u_\tau^0\right\|^{\frac{2 n}{n-2}}+\delta\left\|\nabla u_\tau^0\right\|^4\right) \\
& +C e^{-\sigma_1 t} \int_{-\infty}^t e^{\sigma_1 s}\|h(x, s)\|^2 d s+C
\end{aligned}
\end{equation}
for any $\delta \geq 0$ and $\tau \leq t \in \mathbb{R}$, where $C\neq C(t,\tau,\delta)$ is a positive constant.
\end{Lemma}
$\mathbf{Proof.}$ Taking inner product between $\partial_t u + \rho u$ and $\eqref{1.1-6}_1$ in $L^2(\Omega)$, we obtain
\begin{align}
\frac{1}{2} \frac{d}{d t}& \left(\varepsilon(t)\|\partial_t u+\rho u\|^2-\rho^2 \varepsilon(t)\|u\|^2+\|\nabla u\|^2+\frac{\delta}{2}\|\nabla u\|^4+\rho\|\nabla u\|^2+\lambda\|u\|^2-2\left(G(u), 1\right)\right)\non \\
& -\frac{1}{2} \varepsilon^{\prime}(t)\|\partial_t u\|^2-\rho \varepsilon^{\prime}(t)(\partial_t u, u)-\rho \varepsilon(t)\|\partial_t u\|^2+\rho \|\nabla u\|^2+\delta \rho\| \nabla u\|^4 \non \\
& +\| \nabla \partial_t u\|^2+\lambda \rho\| u \|^2-\rho(g(u), u)  =(h(x, t), \partial_t u+\rho u),\label{3.2-6}
\end{align}
where $\rho >0$ is a constant.

Then from \eqref{3.2-6}, we arrive at
\begin{align}\label{3.3-6}
& \frac{d}{d t} E(t)-\varepsilon^{\prime}(t)\|\partial_t u\|^2-2 \rho \varepsilon^{\prime}(t)(\partial_t u, u)-2 \rho \varepsilon(t) \| \partial_t u\|^2+2 \rho\| \nabla u\|^2\non \\
& +2 \delta \rho\| \nabla u\|^4+2\| \nabla \partial_t u \|^2 +2 \lambda \rho\|u\|^2-2 \rho(g(u), u)+4 \rho^2 \varepsilon(t)(\partial_t u, u)\non \\
& =2\left(h(x, t), \partial_t u+\rho u)+4 \rho^2 \varepsilon(t)(\partial_t u, u), \right.
\end{align}
where
\begin{equation}\label{3.4-6}
\begin{aligned}
E(t)&=\varepsilon(t)\|\partial_t u+\rho u\|^2-\rho^2 \varepsilon(t)\|u\|^2+\|\nabla u\|^2+\frac{\delta}{2}\|\nabla u\|^4+\rho\|\nabla u\|^2\\
&+\lambda\|u\|^2-2(G(u), 1)+2 c_0
\end{aligned}
\end{equation}
and $ c_0>0$ is a constant.

Thanks to \eqref{1.2-6}, \eqref{1.3-6}, the Poincar\'{e} inequality and noticing that $\varepsilon(t)$ is decreasing, then if $\rho$ satisfies $\sqrt{2\lambda} \leq \rho \leq \min \left\{\frac{\lambda_1}{4 L}, \frac{\sqrt{\left(\lambda_1+4 \lambda\right) L}}{2 L}\right\}$, it follows that
\begin{align}\label{3.5-6}
& -\varepsilon^{\prime}(t)\|\partial_t u\|^2-2 \rho \varepsilon^{\prime}(t)(\partial_t u, u)-2 \rho \varepsilon(t)\|\partial_tu\|^2+\frac{\rho}{2}\|\nabla u\|^2 +\|\nabla \partial_tu\|^2\non\\
&+2 \lambda \rho\|u\|^2+4 \rho^2 \varepsilon(t)\left(\partial_t u, u\right) \non\\
& \geq 2 \rho \varepsilon(t)\|\partial_t u+\rho u\|^2.
\end{align}

Furthermore, by the Young inequality, we conclude
\begin{equation}\label{3.6-6}
(h(x, t), \partial_tu+\rho u) \leq \frac{1}{2 \rho}\|h(x, t)\|^2+\frac{\rho}{2}\|\partial_t u+\rho u\|^2
\end{equation}
and
\begin{equation}\label{3.7-6}
2 \rho^2(\partial_t u, u) \leq \frac{4 \rho^2}{\lambda_1}\|\partial_t u\|^2+\frac{\rho^2 \lambda_1}{4}\|u\|^2.
\end{equation}

Inserting $(\ref{3.5-6})-(\ref{3.7-6})$ into $(\ref{3.3-6})$, we deduce
\begin{align}\label{3.8-6}
\frac{d}{dt} & E(t)+\frac{1}{2}\|\nabla \partial_tu\|^2+\frac{\rho}{2}\|\nabla u\|^2+2 \delta \rho\|\nabla u\|^4-2 \rho(g(u), u)+\rho(2 \varepsilon(t)-\rho)\|\partial_t u+\rho u\|^2 \non\\
& +\frac{1}{2}\|\nabla \partial_tu\|^2+\rho{\|\nabla u\|^2}-\frac{8 \rho^2 \varepsilon(t)}{\lambda_1}\|\partial_t u\|^2-\frac{\rho^2 \lambda_1 \varepsilon(t)}{2}\|u\|^2 \non\\
& \leq \frac{1}{\rho}\|h(x, t)\|^2,
\end{align}
where $\chi>0$ is a constant.

Suppose
\begin{equation}\label{3.9-6}
I(t)=\frac{\rho}{2}\|\nabla u\|^2+2 \delta \rho\|\nabla u\|^4-2 \rho(g(u), u)+\rho(2 \varepsilon(t)-\rho)\|\partial_tu+\rho u\|^2-\chi E(t)
\end{equation}
and
\begin{equation}\label{3.10-6}
K(t)=\frac{1}{2}\|\nabla \partial_t u\|^2+\rho\|\nabla u\|^2-\frac{8 \rho^2 \varepsilon(t)}{\lambda_1}\|\partial_t u\|^2-\frac{\rho^2 \lambda_1 \varepsilon(t)}{2}\|u\|^2,
\end{equation}
then we derive
\begin{equation}\label{3.11-6}
\frac{d}{d t} E(t)+\frac{1}{2}\|\nabla \partial_t u\|^2+I(t)+K(t)+ \chi E(t) \leq \frac{1}{\rho}\|h(x, t)\|^2.
\end{equation}

Using \eqref{3.10-6} and the Poincar\'{e} inequality, then if $\rho$ further satisfies $\rho \leq \min \left\{\frac{2}{L}, \frac{\lambda_1 \sqrt{L}}{4 L}\right\}$, we obtain
\begin{equation}\label{3.12-6}
K(t) \geq 0.
\end{equation}

Thanks to \eqref{1.6-6} and \eqref{1.7-6}, we conclude that there exist constants $c_1,c_2,c_3,c_4>0$ such that
\begin{equation}\label{3.13-6}
-(g(u), u) \geq -\gamma(G(u), 1)-c_1\|u\|^2-c_2
\end{equation}
and
\begin{equation}\label{3.14-6}
-({G}(u), 1) \geq -c_3\|u\|^2-c_4.
\end{equation}

Moreover, if $\rho, L, c_3$ satisfy $\lambda_1+\rho \lambda_1-\rho^2 L-2 c_3+\lambda_1 \geq 0$ and let $c_0<c_4$, then by the Poincar\'{e} inequality, we deduce
\begin{equation}\label{3.15-6}
-\rho^2 \varepsilon(t)\|u\|^2+\|\nabla u\|^2+\rho\|\nabla u\|^2+\lambda\|u\|^2-2(G(u), 1)+2 c_0 \geq0,
\end{equation}
which leads to
\begin{equation}\label{3.16-6}
E(t) \geq 0.
\end{equation}

Inserting \eqref{3.4-6} into \eqref{3.9-6}, we obtain
\begin{equation}\label{3.17-6}
\begin{aligned}
I(t)& =\left(\frac{\rho}{2}-\chi-\chi \rho\right)\|\nabla u\|^2+\left(2 \delta \rho-\chi \frac{\delta}{2}\right)\|\nabla u\|^4+(\rho(2 \varepsilon(t)-\rho)-\chi \varepsilon(t))\|\partial_tu+\rho u\|^2 \\
& +\left(\chi \rho^2 \varepsilon(t)-\chi \lambda)\|u\|^2-2 \rho(g(u), u\right)+2 \chi(G(u), 1)-2 \chi c_0.
\end{aligned}
\end{equation}

Thus, from \eqref{3.13-6}, \eqref{3.14-6} and \eqref{3.17-6}, if $\chi$ further fulfills $\chi<\rho \gamma$, we arrive at
\begin{equation}\label{3.18-6}
\begin{aligned}
& -2\rho(g(u),u)+2\chi(G(u),1)\\
& \geq -(2 \rho \gamma c_3-2 \chi c_3+2 \rho c_1)\|u\|^2-(2 \rho \gamma-2 \chi) c_4-2 \rho c_2.
\end{aligned}
\end{equation}

Inserting \eqref{3.18-6} into \eqref{3.17-6}, we conclude
\begin{equation}\label{3.19-6}
\begin{aligned}
I(t)& =\left(\frac{\rho}{2}-\chi-\chi \rho\right)\|\nabla u\|^2+\left(2 \delta \rho-\chi \frac{\delta}{2}\right)\|\nabla u\|^4+\left(2 \rho \varepsilon(t)-\rho^2-\chi \varepsilon(t)\right)\|\partial_t u+\rho u\|^2\\
& +\left(\chi \rho^2 \varepsilon(t)-\chi \lambda-2 \rho \gamma c_3+2 \chi c_3-2 \rho c_1\right) \|u \|^2-(2 \rho \gamma-2 \chi) c_4-2 \rho c_2-2 \chi c_0.
\end{aligned}
\end{equation}

From \eqref{3.19-6}, we derive that if the following inequalities hold
\begin{equation}\label{3.20-6}
\left\{\begin{array}{l}
\frac{\rho}{2}-\chi-\chi \rho \geq 0, \\
2 \delta \rho-\chi \frac{\delta}{2} \geq 0, \\
2 \rho \varepsilon(t)-\rho^2-\chi \varepsilon(t) \geq 0, \\
\chi \rho^2 \varepsilon(t)-\chi \lambda-2 \rho \gamma c_3+2 \chi c_3-2 \rho c_1 \geq 0, \\
-(2 \rho \gamma-2 \chi) c_4-2 \rho c_2-2 \chi c_0 \geq 0,
\end{array}\right.
\end{equation}
then there exists a constant $c_5>0$ such that
\begin{equation}\label{3.21-6}
I(t)\geq -c_5.
\end{equation}

Noticing $c_0<c_4$ and by \eqref{1.2-6}, \eqref{1.3-6} and \eqref{3.20-6}, then \eqref{3.21-6} holds when $c_3 \leq \frac{\lambda}{2}$, $\frac{\sqrt{\left(\lambda-2 c_3\right) L}}{L} \leq \rho \leq 2$ and
\begin{equation}\label{3.22-6}
\max \left\{\frac{2 \rho \gamma c_3+2 \rho c_1}{\rho^2-\lambda+2 c_3}, \frac{2 \rho \gamma c_4+2 \rho c_2}{2 c_4-2 c_0}\right\}<\chi<\min \left\{\frac{\rho}{2(1+\rho)}, 4 \rho, \rho \gamma, 2 \rho-\rho^2\right\}.
\end{equation}

Therefore, from \eqref{3.11-6}, \eqref{3.12-6}, \eqref{3.16-6} and \eqref{3.21-6}, we obtain
\begin{equation}\label{3.23-6}
\frac{d}{d t} E(t) \leq -\chi E(t)+\frac{1}{\rho}\|h(x, t)\|^2+c_5.
\end{equation}

Then taking $0<\sigma_1<\chi$ and by \eqref{3.23-6}, we conclude
\begin{equation}\label{3.24-6}
\frac{d}{d t}\left(e^{\sigma_1 t} E(t)\right) \leq -\left(\chi-\sigma_1\right) e^{\sigma_1 t} E(t)+\frac{e^{\sigma_1 t}}{\rho}\|h(x, t)\|^2+e^{\sigma_1 t} c_5.
\end{equation}

Integrating $(\ref{3.24-6})$ from $\tau$ to $t$, we derive
\begin{equation}\label{3.25-6}
E(t)+(\chi-\sigma_1) e^{-\sigma_1 t} \int_\tau^t e^{\sigma_1 s} E(s) d s \leq e^{-\sigma_1(t-\tau)} E(\tau)+\frac{1}{\rho e^{\sigma_1 t}} \int_\tau^t e^{\sigma_1 s}\|h(x,s)\|^2 d s+\frac{c_5}{\sigma_1}.
\end{equation}

Carrying out some calculations similar to $(\ref{3.7-6})$, we deduce
\begin{equation}\label{3.26-6}
2 \rho \varepsilon(t)|(\partial_t u,u)| \leq \frac{4 \rho}{\lambda_1} \varepsilon(t)\|\partial_t u\|^2+\frac{\rho \lambda_1}{4} \varepsilon(t) \|u\|^2.
\end{equation}

Noticing that $\rho \leq \frac{\lambda_1}{4}$ and by \eqref{1.2-6}, \eqref{1.3-6} and the Poincar\'{e} inequality, we obtain if $L\geq\frac{64}{\lambda_1}$, $\rho$ further satisfies $\rho \geq\frac{\lambda_1 L+\sqrt{\left(\lambda_1^2 L-64 \lambda_1\right) L}}{8 L}$, then there exists a constant
$$
c_6 \geq \max \left\{\frac{\lambda_1}{\lambda_1-4 \rho}, \frac{4 \lambda_1+\rho^2 L-\rho L \lambda_1}{4 \lambda_1}\right\}
$$
such that
\begin{equation}\label{3.27-6}
\|\nabla u\|^2+\varepsilon(t)\left\|\partial_t u\right\|^2 \leq c_6\left(\|\nabla u\|^2+\varepsilon(t)\left\|\partial_t u+\rho u\right\|^2\right).
\end{equation}

Inserting \eqref{3.14-6} into \eqref{3.4-6}, we arrive at
\begin{equation}\label{3.28-6}
E(t) \geq \varepsilon(t)\|\partial_t u+\rho u\|^2-\rho^2 \varepsilon(t)\|u\|^2+\|\nabla u\|^2+\frac{\delta}{2}\|\nabla u\|^4+\rho\|\nabla u\|^2+\lambda\|u\|^2-2 c_3\|u\|^2.
\end{equation}

Besides, inserting \eqref{3.27-6} into \eqref{3.28-6}, we conclude if $c_3 \leq \frac{\lambda-\rho^2 L}{2}$ and $\rho$ further satisfies $\rho \le \frac{{\sqrt {\lambda L} }}{L}$, then
\begin{equation}\label{3.29-6}
c_6^{-1}\left(\|\nabla u\|^2+\varepsilon(t)\|\partial_t u\|^2\right) \leq E(t).
\end{equation}

Similarly, by \eqref{3.26-6}, we obtain
\begin{equation}\label{3.30-6}
\varepsilon(t) \| \partial_t u+\rho u\|^2-\rho^2 \varepsilon(t)\| u \|^2\leq \left(1+\frac{4 \rho}{\lambda_1}\right) \varepsilon(t)\|\partial_t u\|^2+\frac{\rho \lambda_1}{4} \varepsilon(t)\|u\|^2.
\end{equation}

Moreover, from \eqref{1.4-6} and $G(u)=\int_0^u g(s) d s$, we deduce that there exist constants $c_8>c_7>0$ such that
\begin{equation}\label{3.31-6}
\begin{aligned}
-2({G}(u), 1)& \leq 2 c_{7}\left(\int_0^u |u|^{\frac{n+2}{n-2}} d s, 1\right)+2 c_{7}(u, 1) \\
& \leq c_8\|u\|^{\frac{4}{n-2}}.
\end{aligned}
\end{equation}

Then by \eqref{3.4-6}, \eqref{3.30-6} and \eqref{3.31-6}, we derive there exist constants
$$c_9 \geq \max \left\{1+\frac{4 \rho}{\lambda_1}, \frac{\rho L}{4}+\frac{\lambda}{\lambda_1}+1+\rho, \frac{1}{2}, c_{8}\right\}>0$$
and $c_{10}>c_0$ such that
\begin{equation}\label{3.32-6}
E(t) \leq  c_9\left(\|\nabla u\|^2+\varepsilon(t)\|\partial_t u\|^2+\|\nabla u\|^{\frac{2 n}{n-2}}+\delta\|\nabla u\|^4\right)+2 c_{10}.
\end{equation}

Hence, combining \eqref{3.29-6} with \eqref{3.32-6}, we obtain
\begin{equation}\label{3.33-6}
c_6^{-1}\left(\|\nabla u\|^2+\varepsilon(t)\|\partial_t u\|^2) \leq E(t)\right. \leq  c_9\left(\|\nabla u\|^2+\varepsilon(t)\|\partial_t u\|^2+\|\nabla u\|^{\frac{2 n}{n-2}}+\delta\|\nabla u\|^4\right)+2 c_{10}.
\end{equation}

From \eqref{3.16-6}, \eqref{3.25-6} and noticing $0<\sigma_1<\chi$, we deduce
\begin{equation}\label{3.34-6}
c_6 E(t) \leq c_6 e^{-\sigma_1(t-\tau)} E(\tau)+\frac{c_6}{\rho e^{\sigma_1 t}} \int_\tau^t e^{\sigma_1 s}\left\|h(x,s)\right\|^2 d s+\frac{c_5c_6}{\sigma_1}.
\end{equation}

Using \eqref{1.2-6}, \eqref{1.3-6} and \eqref{3.33-6}, we derive
\begin{equation}\label{3.35-6}
E(\tau)\leq c_9\left(\left\|\nabla u_\tau^0\right\|^2+\varepsilon(\tau)\left\| u_\tau^1\right\|^2+\left\|\nabla u_\tau^0\right\|^{\frac{2n}{n-2}}+\delta\left\|\nabla u_\tau^0\right\|^4\right)+2 c_{10}.
\end{equation}

Then inserting \eqref{3.33-6} and \eqref{3.35-6} into \eqref{3.34-6}, we conclude
\begin{equation}\label{3.36-6}
\begin{aligned}
\|\nabla u\|^2+\varepsilon(t)\|\partial_t u\|^2& \leq c_6 c_9 e^{-\sigma_1(t-\tau)}\left(\left\|\nabla u_\tau^0\right\|^2+\varepsilon(\tau)\left\|u_\tau^1\right\|^2+\left\|\nabla u_\tau^0\right\|^{\frac{2 n}{n-2}}+\delta\left\|\nabla u_\tau^0\right\|^4\right) \\
& +\frac{c_6}{\rho e^{\sigma_1t}} \int_\tau^t e^{\sigma_1 s}\|h(x, s)\|^2 d s+2 c_6 c_{10} e^{-\sigma_1(t-\tau)}+\frac{c_5c_6}{\sigma_1}.
\end{aligned}
\end{equation}

Let $c_{11}=c_6 c_9$, $c_{12}=\frac{c_6}{\rho}$ and $c_{13}=2 c_6 c_{10}+\frac{c_5 c_6}{\sigma_1}$, then we obtain
\begin{equation}\label{3.37-6}
\begin{aligned}
\|\nabla u\|^2+\varepsilon(t)\|\partial_t u\|^2& \leq c_{11} e^{-\sigma_1(t-\tau)}\left(\left\|\nabla u_\tau^0\right\|^2+\varepsilon(\tau)\left\|u_\tau^1\right\|^2+\left\|\nabla u_\tau^0\right\|^{\frac{2 n}{n-2}}+\delta\left\|\nabla u_\tau^0\right\|^4\right) \\
& +c_{12} e^{-\sigma_1t}\int_\tau^t e^{\sigma_1 s}\|h(x, s)\|^2 d s+c_{13}.
\end{aligned}
\end{equation}

Consequently, \eqref{3.1-6} holds directly. $\hfill$$\Box$

\begin{Lemma}
Under the assumptions of $\varepsilon(t)$, $g(u)$, $h(x,t)$, $\delta$ and $\lambda$ in Section 1, then for any $\delta \geq 0$, the process $\left\{U_\delta(t, \tau)\right\}_{t \geq \tau}$ generated by problem \eqref{1.1-6} has a pullback absorbing family $\mathcal D_\delta=\left\{D_\delta(t)\right\}_{t \in \mathbb{R}}$ that satisfies \eqref{1b-6}.
\label{lem3.3-6}
\end{Lemma}
$\mathbf{Proof.}$ Let
\begin{equation}\label{3.38-6}
B(t)=\left(c_{14} e^{-\sigma_1 t} \int_{-\infty}^{t} e^{\sigma_1 s}\|h(x, s)\|^2 d s+c_{14}\right)^\frac{1}{2},
\end{equation}
where $c_{14} \geq \max \left\{2 c_{12}, 2 c_{13}\right\}$ and $\sigma_1>0$ is the same as in Lemma \ref{lem3.2-6}.

In addition, suppose $\mathcal D_\delta=\left\{D_\delta(t)\right\}_{t \in \mathbb{R}}$ satisfies
\begin{equation}\label{3.39-6}
D_{\delta}(t)=\left\{\left(u^0 ; u^1\right) \in X_t\,\big{|}\left\|\left(u^0 ; u^1\right)\right\|_{X_t} \leq B(t)\right\}.
\end{equation}

Assume $\left(u_{\tilde t-\tau}^0, u_{\tilde t-\tau}^1\right) \in D_{\delta}\left(\tilde t-\tau\right)$ for any $\tilde t \in \mathbb R$, then by $\eqref{3.37-6}-\eqref{3.39-6}$, we deduce there exists a constant $\sigma_2$ satisfying $\sigma_1<\sigma_2<\chi$ such that
\begin{align}\label{3.40-6}
&\left\|U_{\delta}\left(\tilde t, \tilde t-\tau\right) \left(u_{\tilde t-\tau}^0, u_{\tilde t-\tau}^1\right)\right\|_{X_t}^2=\|\nabla u(\tilde t)\|^2+\varepsilon(\tilde t)\|\partial_t u(\tilde t)\|^2\non\\
&\leq c_{11} e^{-\sigma_2 \tau}\left(1+\left(B(\tilde t-\tau)\right)^{\frac{2 n}{n-2}}+\delta\left(B\left(\tilde t-\tau\right)\right)^4 \right) +c_{12} e^{-\sigma_1 \tilde t} \int_{ \tilde t-\tau}^{ \tilde t} e^{\left(\sigma_2-\sigma_1\right) s} e^{\sigma_1 s}\|h(x, s)\|^2 d s+c_{13}\non \\
& \leq c_{14} e^{-\sigma_2 \tau}\left(\left(B(\tilde t-\tau)\right)^{\frac{2 n}{n-2}}+\delta\left(B\left(\tilde t-\tau\right)\right)^4 \right) +c_{15} e^{-\sigma_1 \tilde t} \int_{ -\infty}^{ \tilde t}  e^{\sigma_1 s}\|h(x, s)\|^2 d s+c_{13},
\end{align}
where $c_{14} \geq 2 c_{11}$ and $c_{15}\geq 2 c_{12}$ are constants.

Hence, it follows from \eqref{3.40-6} that for any $\tilde t \in \mathbb R$ there exists a $T=T(\tilde t)>0$ such that
\begin{equation}\label{3.41-6}
U_\delta(\tilde{t}, \tilde{t}-\tau) D_\delta(\tilde{t}-\tau) \subset D_\delta(\tilde{t})
\end{equation}
for any $\tau>T$ and any $\delta \geq 0$. $\hfill$$\Box$

\begin{Lemma}\label{lem3.4-6}
Under the assumptions of $\varepsilon(t)$, $g(u)$, $h(x,t)$, $\delta$ and $\lambda$ in Section 1, assume $\left\|\nabla u_{\tau}^0\right\|^2+\varepsilon(\tau)\|u_{\tau}^1\|^2 \leq \widetilde C$, then the following inequality hold
\begin{equation}\label{3.42-6}
\begin{aligned}
& \varepsilon(t)\left\|\partial_{tt} u(t)\right\|_{-1}^2+\left\|\nabla \partial_t u(t)\right\|^2+\int_t^{t+1}\left\|\partial_{tt}  u(s)\right\|^2 d s \\
& \leq Z\left(\frac{1}{(t-\tau)^2} \int_{\tau}^T(s-\tau)^2\left\|\partial_t h(x, s)\right\|^2 d s+\int_{\tau}^T\left\|\partial_t h(x, s)\right\|^2 d s\right),
\end{aligned}
\end{equation}
where $t \in (\tau, T)$ and $Z(\cdot)=Z(L, \lambda, \rho, \widetilde C)>0$.
\end{Lemma}
$\mathbf{Proof.}$ Integrating \eqref{3.23-6} from $\tau$ to $t$, we conclude
\begin{equation}\label{3.43-6}
\max _{t \in[\tau, T]}\|\nabla u(t)\|^2+\int_{\tau}^T\left\|\nabla \partial_t u(t)\right\|^2 d t \leq c_{16},
\end{equation}
where $c_{16}>0$ is a constant.

Taking the derivative of $\eqref{1.1-6}_1$ with respect to $t$ and assuming $w=\partial_t u$, we obtain
\begin{equation}\label{3.44-6}
\begin{aligned}
& \varepsilon^{\prime}(t) \partial_tw+\varepsilon(t) \partial_{tt} w-2 \delta(\nabla u, \nabla w) \Delta u-\left(1+\delta \|\nabla u\|^2\right) \Delta w-\Delta \partial_tw+\lambda w \\
& =g^{\prime}(u) w+\partial_t h(x, t).
\end{aligned}
\end{equation}

Choosing $(-\Delta)^{-1} \partial_t w+\rho w$ as a test function of \eqref{3.44-6}, we derive
\begin{align}\label{3.45-6}
& \frac{1}{2} \frac{d}{d t}\left(\varepsilon(t)\left\|(-\Delta)^{-\frac{1}{2}} \partial_t w\right\|^2+2 \rho \varepsilon(t)\left(\partial_t w, w\right)+\|w\|^2+\rho\|\nabla w\|^2+\lambda\left\|(-\Delta)^{-\frac{1}{2}} w\right\|^2\right) \non\\
& +\frac{1}{2} \varepsilon^{\prime}(t)\left\|(-\Delta)^{-\frac{1}{2}} \partial_t w\right\|^2-\rho \varepsilon(t)\|\partial_t w\|^2+2 \delta(\nabla u, \nabla w)\left(-\Delta u,(-\Delta)^{-1} \partial_t w\right) +\lambda \rho\|w\|^2 \non\\
& +2 \delta \rho(\nabla u, \nabla w)^2+\delta\|\nabla u\|^2\left(-\Delta w,(-\Delta)^{-1} \partial_t w\right)+\rho\left(1+\delta\|\nabla u\|^2\right)\|\nabla w\|^2+\left\|\partial_t w\right\|^2\non\\
& =\left(g^{\prime}(u) w,(-\Delta)^{-1} \partial_t w\right)+\rho\left(g^{\prime}(u) w, w\right)+\left(\partial_t h(x, t),(-\Delta)^{-1} \partial_t w+\rho w\right).
\end{align}

Thanks to the Cauchy and Young inequalities, we deduce
\begin{align}\label{3.46-6}
&2 \delta(\nabla u, \nabla w)\left(-\Delta u,(-\Delta)^{-1} \partial_t w\right)\non\\
&= \delta\|\nabla u\|^2\left(2(\nabla w, (-\Delta)^{-\frac{1}{2}} \partial_t w)\right) \non\\
&\leq \delta\|\nabla u\|^2\left(\|\nabla w\|^2+\left\|(-\Delta)^{-\frac{1}{2}} \partial_t w\right\|^2\right)
\end{align}
and
\begin{align}\label{3.47-6}
&\delta\|\nabla u\|^2\left(-\Delta w,(-\Delta)^{-1} \partial_t w\right)\non\\
&= \delta \|\nabla u\|^2 \left(\nabla w, (-\Delta)^{-\frac{1}{2}} \partial_t w\right) \non\\
&\leq \frac{\delta}{2}\|\nabla u\|^2\left(\|\nabla w\|^2+\left\|(-\Delta)^{-\frac{1}{2}} \partial_t w\right\|^2\right).
\end{align}

Besides, from \eqref{1.4-6} and the Young inequality, we conclude there exists a constant $c_{17}>0$ such that
\begin{equation}\label{3.48-6}
\begin{aligned}
\left|\left(g^{\prime}(u) w,(-\Delta)^{-1} \partial_t w\right)\right| \leq c_{17}\left(1+\|\nabla u\|^{\frac{4}{n-2}}\right)\left(\|\nabla w\|^2+\left\|(-\Delta)^{-\frac{1}{2}} \partial_t w\right\|^2\right).
\end{aligned}
\end{equation}

Furthermore, it follows from \eqref{1.5-6} that
\begin{equation}\label{3.49-6}
-\rho(g^{\prime}(u) w, w) \leq \rho k\|w\|^2.
\end{equation}

Thanks to the Young inequality, we derive
\begin{equation}\label{3.50-6}
\left(\partial_t h(x, t),(-\Delta)^{-1} \partial_t w+\rho w\right) \leq \frac{1}{2} \| \partial_t h(x,t)\|^2+\left\|(-\Delta)^{-1} \partial_t w\right\|^2+\rho^2\|w\|^2.
\end{equation}

Let
\begin{equation}\label{3.51-6}
L(t)=\varepsilon(t)\left\|(-\Delta)^{-\frac{1}{2}} \partial_t w\right\|^2+2 \rho \varepsilon(t)\left(\partial_t w, w\right)+\|w\|^2+\rho\|\nabla w\|^2+\lambda\left\|(-\Delta)^{-\frac{1}{2}} w\right\|^2,
\end{equation}
then from Lemmas \ref{lem3.2-6} and \ref{lem3.3-6}, we deduce there exist constants $0<c_{18}\leq \min\{1,\rho\}$ and $c_{19}\geq \max\{1,\rho\}$ depending on $L, \lambda$ and $\rho$ such that
\begin{equation}\label{3.52-6}
c_{18}\left(\varepsilon(t)\left\|(-\Delta)^{-\frac{1}{2}} \partial_t w\right\|^2+\|\nabla w\|^2\right) \leq L(t) \leq c_{19}\left(\varepsilon(t)\left\|(-\Delta)^{-\frac{1}{2}} \partial_t w\right\|^2+\|\nabla w\|^2\right).
\end{equation}

Inserting \eqref{3.51-6} into \eqref{3.45-6}, then by \eqref{1.2-6}, \eqref{1.3-6} and $\eqref{3.46-6}-\eqref{3.50-6}$, we arrive at
\begin{align}\label{3.53-6}
\frac{d}{d t} L(t)+2\|\partial_tw\|^2 &\leq (L+2)\left\|(-\Delta)^{-\frac{1}{2}} \partial_tw\right\|^2+2\rho L\|\partial_t w\|^2\non\\
&+2c_{17}\left(1+\|\nabla u\|^{\frac{4}{n-2}}\right)\left(\|\nabla w\|^2+\left\|(-\Delta)^{-\frac{1}{2}} \partial_t w\right\|^2\right)\non\\
&+\|\partial_t h(x, t) \|^2.
\end{align}

Moreover, from Lemma \ref{lem3.2-6} and noticing $n\geq 3$, we obtain that there exists a constant $c_{20}>0$ such that
\begin{equation}\label{3.54-6}
\begin{aligned}
\frac{d}{d t} L(t)+ \|\partial_tw\|^2 &\leq (L+2)\left\|(-\Delta)^{-\frac{1}{2}} \partial_tw\right\|^2+2\rho L\|\partial_t w\|^2\\
&+c_{20}\left(\|\nabla w\|^2+\left\|(-\Delta)^{-\frac{1}{2}} \partial_t w\right\|^2\right)+ \|\partial_t h(x, t) \|^2.
\end{aligned}
\end{equation}

Hence, there exists a constant $c_{21}>\max \left\{L+2+c_{20}, 2 \rho L\right\}$ depending on $\widetilde{C}, L, \rho$ and $\lambda$ such that \begin{equation}\label{3.55-6}
\begin{aligned}
\frac{d}{d t} L(t)+ \|\partial_tw\|^2 &\leq c_{21}L(t)+ \|\partial_t h(x, t) \|^2.
\end{aligned}
\end{equation}

Noticing $w=\partial_t u$ and inserting \eqref{3.51-6} into \eqref{3.55-6}, we obtain
\begin{equation}\label{3.56-6}
\begin{aligned}
& \frac{d}{d t}L(t) \leq c_{21} L(t)+\|\partial_t h(x, t)\|^2.
\end{aligned}
\end{equation}

By the Gronwall inequality, we deduce that there exist constants $c_{22}>c_{23}>0$ depending on $\widetilde C,L,\rho$ and $\lambda$ such that
\begin{equation}\label{3.57-6}
L(t) \leq c_{22} L(\tau)+c_{23} \int_\tau^T\|\partial_t h(x, s)\|^2 d s.
\end{equation}

Therefore, there exists a constant $c_{24} \geq \max\{c_{22},c_{23}\}$ depending on $\widetilde C,L,\rho$ and $\lambda$ such that
\begin{equation}\label{3.58-6}
\begin{aligned}
& \varepsilon(t)\left\|\partial_{tt} u(t)\right\|_{-1}^2+\left\|\nabla \partial_t u(t)\right\|^2+\int_t^{t+1}\left\|\partial_{tt}  u(s)\right\|^2 d s \\
& \leq  c_{24} \left(L(\tau)+c_{23} \int_\tau^T\|\partial_t h(x, s)\|^2 d s\right).
\end{aligned}
\end{equation}

Then combining \eqref{3.51-6} with \eqref{3.58-6}, it follows that there exists a constant $c_{25}>0$ such that
\begin{equation}\label{3.59-6}
 L(t)\leq c_{25}\left(\varepsilon(t)\left\|\partial_{tt} u(t)\right\|_{-1}^2+\left\|\nabla \partial_t u(t)\right\|^2\right).
\end{equation}

Thus, thanks to \eqref{3.58-6} with \eqref{3.59-6}, we conclude there exists a function $Z(\cdot)=Z(L, \lambda, \rho, \widetilde C)>c_{24}c_{25}>0$ such that for any $t \in [\tau, T-1]$,
\begin{equation}\label{3.60-6}
\begin{aligned}
& \varepsilon(t)\left\|\partial_{tt} u(t)\right\|_{-1}^2+\left\|\nabla \partial_t u(t)\right\|^2+\int_t^{t+1}\left\|\partial_{tt}  u(s)\right\|^2 d s \\
& \leq Z\left(\varepsilon(\tau)\left\|\partial_{tt} u(\tau)\right\|_{-1}^2+\left\|\nabla \partial_t u(\tau)\right\|^2+\int_{\tau}^T\left\|\partial_t h(x, s)\right\|^2 d s\right).
\end{aligned}
\end{equation}

From \eqref{3.52-6} and \eqref{3.56-6}, we derive
\begin{equation}\label{3.61-6}
\begin{aligned}
\frac{d}{d t}\left((t-\tau)^2 L(t)\right)+(t-\tau)^2\|\partial_tw\|^2 &\leq 2(t-\tau) c_{19}\left(\varepsilon(t)\left\|(-\Delta)^{-\frac{1}{2}} \partial_t w\right\|^2+\|\nabla w\|^2\right)\\
& +c_{21}(t-\tau)^2L(t)+(t-\tau)^2\|\partial_th(x,t)\|^2.
\end{aligned}
\end{equation}

By \eqref{3.51-6} and the embedding theorem, we obtain that there exist constants $c_{26},c_{27}>0$ such that
\begin{equation}\label{3.62-6}
\begin{aligned}
& 2(t-\tau) c_{19}\left(\varepsilon(t)\left\|(-\Delta)^{-\frac{1}{2}} \partial_t w\right\|^2+\|\nabla w\|^2\right) \\
 &\leq c_{26}\left(1+(t-\tau)^2\|\nabla \partial_t u\|^2 L(t)\right)+\frac{1}{2}(t-\tau)^2\|\partial_t w\|^2+c_{27}.
\end{aligned}
\end{equation}

Then inserting \eqref{3.62-6} into \eqref{3.61-6}, it follows that there exist constants $c_{28},c_{29}>0$ such that
\begin{equation}\label{3.63-6}
\frac{d}{d t}\left((t-\tau)^2 L(t)\right) \leq c_{28}(t-\tau)^2 L(t)+c_{29}(t-\tau)^2\|\partial_t h(x, t)\|^2.
\end{equation}


Using the Gronwall inequality, we deduce that there exists a constant $C>0$ such that
\begin{equation}
\begin{aligned}
(t-\tau)^2 L(t) & \leq c_{29} \int_\tau^t e^{\int_\tau^t c_{28} d r}(s-\tau)^2\|\partial_t h(x, s)\|^2 d s \\
& \leq c_{29}C \int_\tau^T(s-\tau)^2\|\partial_t h(x, s)\|^2 d s.
\end{aligned}
\end{equation}

Hence, it follows that
\begin{equation}\label{3.63.2-6}
L(t) \leq \left.\frac{c_{29}C}{(t-\tau)^2} \int_{\tau}^T(s-\tau)^2\left\|\partial_t h(x, s)\right\|^2\right. d s.
\end{equation}

Inserting \eqref{3.63.2-6} into \eqref{3.60-6}, then \eqref{3.42-6} follows directly. $\hfill$$\Box$

Next, we establish in Lemma \ref{lem3.5-6} the existence of a pullback absorbing family in $X_t$ for problem \eqref{1.1-6}, which is helpful to derive the existence of pullback attractors.
\begin{Lemma} \label{lem3.5-6}
Under the assumptions of Lemma $\ref{lem3.4-6}$, if there exists a $\tau_{\delta}>0$ such that
\begin{equation}\label{3.64-6}
U_{\delta}(t, t-\tau) D_{\delta}(t-\tau) \subset D_{\delta}(t)
\end{equation}
and
\begin{equation}\label{3.65-6}
U_{\delta}(t-1, t-\tau) D_{\delta}(t-\tau) \subset D_{\delta}(t-1)
\end{equation}
for any $\tau\geq\tau_{\delta}$ and $t \in \mathbb R$, and assume $\widetilde {\mathcal{D}}_{\delta}=\{\widetilde {\mathcal{D}}_{\delta}(t)\}_{t \in \mathbb{R}}$ satisfies
\begin{equation}\label{3.66-6}
\widetilde {\mathcal{D}}_{\delta}(t)= \mathop  {\overline {\mathop  \bigcup \limits_{\tau  \ge \tau_{\delta}} U_{\delta}(t,t - \tau )D_{\delta}(t - \tau )} ^{X_t}},
\end{equation}
then the following properties hold

\noindent(i) $\widetilde {\mathcal{D}}_{\delta}=\{\widetilde {\mathcal{D}}_{\delta}(t)\}_{t \in \mathbb R}$ is a pullback absorbing family in $X_t$ and satisfies
\begin{equation}\label{3.67-6}
U_{\delta}(t, t-\tau) \widetilde {D}_{\delta}(t-\tau) \subset \widetilde {D}_{\delta}(t), \quad \forall \, \tau \geq T;
\end{equation}

\noindent(ii) $\widetilde {D}_{\delta}(t)$ is bounded in $H_0^{1}(\Omega)\times H_0^{1}(\Omega)$.
\end{Lemma}
$\mathbf{Proof.}$ Suppose $B_1 \subset X_t$ is bounded and choose a time $\tau_1>0$ large enough, then noticing $\mathcal{D}_{\delta}=\{D_{\delta}(t)\}_{t \in \mathbb{R}}$ is a pullback absorbing family in $X_t$, we derive
\begin{equation}\label{3.68-6}
U_{\delta}\left(t-\tau_{\delta}-1, t-\tau\right) B_1 \subseteq D_{\delta}\left(t-\tau_{\delta}-1\right), \quad \forall\, \tau \geq \tau_1.
\end{equation}

By \eqref{3.66-6} and \eqref{3.68-6}, then we conclude for any $\tau\geq\tau_{\delta}$,
\begin{align}\label{3.69-6}
U_\delta(t, t-\tau) B_1 & =U_\delta\left(t, t-\tau_\delta-1\right) U_\delta\left(t-\tau_\delta-1, t-\tau\right) B_1 \non\\
& \subseteq U_\delta\left(t, t-\tau_\delta-1\right) D_\delta\left(t-\tau_\delta-1\right)\non \\
&\subseteq \widetilde {D}_\delta(t).
\end{align}

Similarly, we deduce
\begin{equation}\label{3.70-6}
\widetilde {D}_\delta(t) \subset {D}_\delta(t), \quad \forall \, t\in \mathbb R.
\end{equation}

Therefore, it follows that
\begin{equation}\label{3.71-6}
U_\delta(t, t-\tau) \widetilde{D}_\delta(t-\tau) \subseteq U_\delta(t, t-\tau) D_\delta(t-\tau) \subseteq \widetilde{D}_\delta(t), \quad \forall \,\tau \geq\tau_{\delta}.
\end{equation}

Then from $\eqref{3.64-6}-\eqref{3.66-6}$, we obtain
\begin{align}\label{3.72-6}
\widetilde{D}_{\delta}(t) &\subseteq U_\delta(t, t-\tau) D_\delta(t-\tau)\non\\
&=U_\delta(t, t-1) U_\delta(t-1, t-\tau) D_\delta(t-\tau) \non\\
&\subseteq U_\delta(t, t-1) D_\delta(t-1)
\end{align}
for any $\tau \geq\tau_{\delta}$ and $t\in\mathbb R$.

Consequently, by \eqref{3.43-6}, it follows that $\widetilde{D}_\delta(t)$ satisfies Lemma \ref{lem3.5-6} (ii).$\hfill$$\Box$

Before proving the process $U_{\delta}(\cdot, \cdot)$ generated by problem \eqref{1.1-6} is pullback $\mathcal{D}_{\delta}$-asymptotically compact in the time-dependent space $X_t$, we verify the following lemma in advance.

\begin{Lemma} \label{lem3.6-6}
Under the assumptions of Lemma $\ref{lem3.5-6}$, suppose $\widetilde {\mathcal{D}}_{\delta}=\{\widetilde {\mathcal{D}}_{\delta}(t)\}_{t \in \mathbb R}$ satisfies \eqref{3.66-6} and the process $U_{\delta}(\cdot, \cdot)$ is pullback $\widetilde{\mathcal{D}}_{\delta}$-asymptotically compact in the time-dependent space $X_t$, then $U_{\delta}(\cdot, \cdot)$ is also pullback ${\mathcal D_{\delta}}$-asymptotically compact in $X_t$. In addition, it follows that
\begin{equation}\label{3.73-6}
\Lambda(\widetilde{\mathcal{D}}_{\delta}, t)=\Lambda(\mathcal{D}_{\delta}, t).
\end{equation}
\end{Lemma}
$\mathbf{Proof.}$ Let $i \in \mathbb N$ and assume the sequences $r_i \rightarrow +\infty, \tau_i \rightarrow +\infty$ are increasing  and $x_i \in D_{\delta}\left(t-\tau_i\right)$.

Then from \eqref{3.66-6}, we conclude
\begin{equation}\label{3.74-6}
\widetilde {\mathcal{D}}_{\delta}(t-r_i)= \mathop {\overline {\mathop \bigcup \limits_{\tau  \ge \tau_{\delta,t-r_i}} U_{\delta}(t-r_i,t - \tau )D_{\delta}(t - \tau )} ^{X_t}}.
\end{equation}

Besides, let $j \in \mathbb N$, and select a sequence $\tau_{i_j} \in\{\tau_i\}_{i \in \mathbb{N}}$ satisfies $\tau_{i_j} \geq \tau_{\delta, t-r_j,} \tau_{i_j} \geq  \tau_{\tau_{i_{j-1}}}$ and $\tau_{i_j} \rightarrow+\infty$ for any $r_j \in\left\{r_i\right\}_{i \in \mathbb N}$.

Suppose $a_j=U_\delta\left(t-r_j, t-\tau_{i_j}\right) x_{i_j}$, then by $\left.x_i \in D_\delta ( t-\tau_i\right)$ and \eqref{3.66-6}, we derive
\begin{equation}\label{3.75-6}
a_j\in U_{\delta}\left(t-r_j, t-\tau_{i_j}\right) D_\delta\left(t-\tau_{i_j}\right) \subseteq \widetilde{D}_{\delta}\left(t-r_j\right)
\end{equation}
and
\begin{equation}\label{3.76-6}
\begin{aligned}
U_\delta\left(t, t-\tau_{i_j}\right) x_{i_j} & =U_\delta\left(t, t-r_j\right) U_\delta\left(t-r_j, t-\tau_{i_j}\right) x_{i_j} \\
& =U_\delta\left(t, t-r_j\right) a_j.
\end{aligned}
\end{equation}

From \eqref{3.75-6} and \eqref{3.76-6}, we obtain
\begin{equation}\label{3.77-6}
U_\delta\left(t, t-\tau_{i_j}\right) x_{i_j}=U_\delta\left(t, t-r_j\right) a_j \subseteq U_\delta\left(t, t-r_j\right) \widetilde{D}_{\delta}\left(t-r_j\right).
\end{equation}

Noticing the process $U_{\delta}(\cdot, \cdot)$ is pullback $\widetilde{\mathcal{D}}_{\delta}$-asymptotically compact in the time-dependent space $X_t$, it follows that $\left\{U_\delta\left(t, t-\tau_{i_j}\right) x_{i_j}\right\}$ is precompact in $X_t$.

Therefore, by Definition \ref{def2.3-6}, we conclude that $U_{\delta}(\cdot, \cdot)$ is pullback ${\mathcal D_{\delta}}$-asymptotically compact in $X_t$.

Furthermore, thanks to \eqref{3.66-6}, we deduce
\begin{equation}\label{3.78-6}
\widetilde {D}_\delta(t) \subset {D}_\delta(t),\quad\forall\, t\in\mathbb R,
\end{equation}
then from Lemma \ref{lem2.5-6}, we arrive at
\begin{equation}\label{3.79-6}
 \Lambda(\widetilde{\mathcal D}_{\delta}, t) \subset \Lambda\left(\mathcal D_\delta, t\right),\quad \forall \,t \in \mathbb{R}.
\end{equation}

Suppose $a \in \Lambda({D}_{\delta}, t)$, then there exist sequences $\tau_i \rightarrow+\infty$ and  $x_i \in D_{\delta}\left(t-\tau_i\right)$ such that
\begin{equation}\label{3.80-6}
U_\delta\left(t, t-\tau_i\right) x_i \in U_\delta\left(t, t-\tau_i\right) D_\delta\left(t-\tau_i\right) \rightarrow a \quad\text { as } i \rightarrow+\infty.
\end{equation}

Similarly, we can choose sequences $r_j \rightarrow + \infty$ and $\tau_{i_j} \in\left\{\tau_i\right\}_{i \in \mathbb N}$ such that
\begin{equation}\label{3.81-6}
a_j=U_{\delta}\left(t-r_j, t-\tau_{i_j}\right) x_{i_j} \subset \widetilde{D}_{\delta}\left(t-r_j\right).
\end{equation}

Then by \eqref{3.80-6}, we derive
\begin{align}\label{3.82-6}
 U_\delta\left(t, t-r_j\right) a_j &=U_\delta\left(t, t-r_j\right) U_\delta\left(t-r_j, t-\tau_{i_j}\right) x_{i_j} \non\\
& =U_\delta\left(t, t-\tau_{i_j}\right) x_{i_j} \rightarrow a.
\end{align}

Furthermore, from \eqref{3.74-6} and \eqref{3.75-6}, we obtain
\begin{equation}\label{3.83-6}
U_\delta\left(t, t-r_j\right) a_j \in U_\delta\left(t, t-r_j\right) \widetilde{D}_{\delta}\left(t-r_j\right) \subseteq \Lambda(\widetilde{\mathcal D}_{\delta}, t).
\end{equation}

Hence, thanks to \eqref{3.81-6} and \eqref{3.82-6}, we conclude
\begin{equation}\label{3.84-6}
U_\delta\left(t, t-r_j\right) a_j \rightarrow a \in \Lambda(\widetilde{\mathcal D}_{\delta}, t).
\end{equation}

Consequently, \eqref{3.73-6} holds immediately.
$\hfill$$\Box$

\begin{Theorem}\label{th3.7-6}
Under the assumptions of Lemma $\ref{lem3.6-6}$, then for any $\delta \geq 0$, there exists a pullback attractor $\widetilde{\mathcal A}_\delta=\{\widetilde{A}(t)\}_{t \in \mathbb{R}}$ for the process $U_{\delta}(\cdot,\cdot)$ generated by problem \eqref{1.1-6} in the time-dependent space $X_t$, which satisfies
\begin{equation}\label{3.85-6}
\widetilde{\mathcal A}_{\delta}(t)=\Lambda(\widetilde{\mathcal D}_{\delta}, t), \quad \forall \, t\in\mathbb R.
\end{equation}
\end{Theorem}
$\mathbf{Proof.}$ From Lemmas \ref{lem2.5-6} and \ref{lem3.3-6}, we only need to prove the process $U_{\delta}(\cdot,\cdot)$ is pullback $\mathcal D_{\delta}$-asymptotically compact in $X_t$.

Let $\delta \in [0,1]$, then by Lemma \ref{lem3.3-6} and \eqref{3.43-6}, we deduce
\begin{equation}\label{3.86-6}
\max _{t \in[t_0-\tau, t_0]}\|\nabla u(t)\|^2+\int_{t_0-\tau}^t\left\|\nabla \partial_t u(t)\right\|^2 d t \leq c_{16}
\end{equation}
for any $\left(u_{t_0-\tau}^0, u_{t_0-\tau}^1\right) \in \widetilde{D}_\delta\left(t_0-\tau\right)$, where the constant $c_{16}$ is the same as in \eqref{3.43-6}.

Decompose the solutions to problem \eqref{1.1-6} as
\begin{align}\label{3.87-6}
& U_\delta\left(t, t_0-\tau\right)\left(u_{t_0-\tau}^{0} ; u_{t_0-\tau}^1\right) \non\\
&= \left(U_{a, \delta}\left(t, t_0-\tau\right) u_{t_0-\tau}^0 ; U_{b, \delta}\left(t, t_0-\tau\right) u_{t_0-\tau}^1\right)\non \\
&= \left(u(t) ; \partial_t u(t)\right),
\end{align}
where
\begin{equation}\label{3.88-6}
U_{a, \delta}\left(t, t_0-\tau\right) u_{t_0-\tau}^{0}=U_{a_1,\delta}\left(t, t_0-\tau\right) u_{t_0-\tau}^0+U_{a_2, \delta}\left(t, t_0-\tau\right) u_{t_0-\tau}^0
\end{equation}
with
\begin{equation}\label{3.89-6}
U_{a_1, \delta}\left(t_1, t_0-\tau\right) u_{t_0-\tau}^0=u_1(t)
\end{equation}
and
\begin{equation}\label{3.90-6}
U_{a_2, \delta}\left(t, t_0-\tau\right) u_{t_0-\tau}^0=u_2(t).
\end{equation}

Then $U_{a_1, \delta}\left(t, t_0-\tau\right) u_{t_0-\tau}^0$ and $U_{a_2, \delta}\left(t, t_0-\tau\right) u_{t_0-\tau}^0$ satisfy the following equations
\begin{equation}\label{3.91-6}
\begin{cases}-\left(1+\delta \| \nabla u\|^2\right) \Delta u_1-\Delta \partial_t u_1+\lambda u_1=\varphi\left(u_1+u_2\right)-\varphi\left(u_2\right) & \text { in } \Omega \times(\tau,+\infty), \\ u_1(x, t)=0 & \text { on } \partial \Omega \times(\tau,+\infty), \\ u_1(x, t-\tau)=u_{t-\tau}^0, & \,\, x \in \Omega,\end{cases}
\end{equation}
and
\begin{equation}\label{3.92-6}
\begin{cases}-\left(1+\delta\|\nabla u\|^2\right) \Delta u_2-\Delta \partial_t u_2+\lambda u_2=\varphi\left(u_2\right)+\psi & \text { in } \Omega \times(\tau,+\infty), \\ u_2(x, t)=0 & \text { on } \partial \Omega \times(\tau,+\infty), \\ u_2(x, t-\tau)=0, &\,\, x\in\Omega,\end{cases}
\end{equation}
respectively, where $\varphi(s)=g(s)-k s$ and $\varphi^{\prime}(s)\leq0$ for any $s\in\mathbb R$,
$\psi=-\varepsilon(t) \partial_{tt} u+k u+h(x, t)$, and $k$ is the same as in \eqref{1.5-6}.

Choosing $u_1$ as a test function of $\eqref{3.91-6}_1$, we derive
\begin{equation}\label{3.93-6}
\frac{d}{d t}\left\|\nabla u_1\right\|^2+2\left\|\nabla u_1\right\|^2+2 \delta\left\|\nabla u \right\|^2 \|\nabla u_1\left\|^2+2 \lambda\right\| u_1 \|^2 \leq 0.
\end{equation}

Then by the Gronwall inequality, \eqref{3.89-6} and $\eqref{3.91-6}_3$, we obtain
\begin{equation}\label{3.94-6}
\left\|\nabla u_1(t)\right\|^2=\left\|\nabla U_{a_1, \delta} \left(t,t_0-\tau\right) u_{t_0-\tau}^0\right\|^2 \leq e^{-2\left(t-t_0+\tau\right)}\left\|\nabla u_{t_0-\tau}^0\right\|^2.
\end{equation}

In addition, choosing $-\Delta u_2$ as a test function of $\eqref{3.92-6}_1$, we conclude
\begin{equation}\label{3.95-6}
\begin{aligned}
 &2\left\|\Delta u_2\right\|^2+2 \delta\|\nabla u\|^2\| \Delta u_2\|^2+\frac{d}{d t}\| \Delta u_2\|^2+2 \lambda\| \nabla u_2 \|^2 \\
&=2\left(\varphi\left(u_2\right),-\Delta u_2\right) +2\left(-\varepsilon(t) \partial_{tt} u+k u+h\left(x, t\right),-\Delta u_2\right).
\end{aligned}
\end{equation}

Thanks to \eqref{1.4-6}, $\varphi(s)=g(s)-k s$ and the Young inequality, we arrive at
\begin{equation}\label{3.96-6}
\begin{aligned}
2|(\varphi(u_2),-\Delta u_2)| &\leq  4\left\|g\left(u_2\right)\right\|^2+\frac{1}{2}\left\|\Delta u_2\right\|^2+4 k^2\left\|u_2\right\|^2 \\
& \leq C\left(1+\left\|u_2\right\|^{\frac{2 n+4}{n-2}}\right)+\frac{1}{2}\left\|\Delta u_2\right\|^2+4 k^2\left\|u_2\right\|^2.
\end{aligned}
\end{equation}

Similarly, by \eqref{1.2-6}, \eqref{1.3-6} and the Young inequality, we deduce
\begin{equation}\label{3.97-6}
\begin{aligned}
&2\left(-\varepsilon(t) \partial_{tt} u+k u+h\left(x, t\right),-\Delta u_2\right)\\
& \leq2 L^2\|\partial_{tt}u\|^2+\frac{3}{2}\left\|\Delta u_2\right\|^2+2 k^2\|u\|^2+2\|h(x, t)\|^2.
\end{aligned}
\end{equation}

Inserting \eqref{3.96-6}, \eqref{3.97-6} into \eqref{3.95-6}, we obtain
\begin{equation}\label{3.98-6}
\begin{aligned}
 &2 \delta\|\nabla u\|^2\| \Delta u_2\|^2+\frac{d}{d t}\| \Delta u_2\|^2+2 \lambda\| \nabla u_2 \|^2 \\
&\leq C\left(1+\left\|u_2\right\|^{\frac{2 n+4}{n-2}}\right)+4 k^2\left\|u_2\right\|^2+2 L^2\|\partial_{tt}u\|^2+2 k^2\|u\|^2+2\|h(x, t)\|^2.
\end{aligned}
\end{equation}

Using Lemma \ref{lem3.3-6} and noticing $u=u_1+u_2$, we derive there exist constants $c_{30},c_{31},c_{32}>0$ depending on $\delta$ and $\lambda$ such that
\begin{equation}\label{3.99-6}
\frac{d}{d t}\left\|\Delta u_2\right\|^2+c_{30}\left\|\Delta u_2\right\|^2 \leq c_{31}+c_{32}\left\|h\left(x,t\right)\right\|^2.
\end{equation}

Let $c_{33}=\max\{c_{31},c_{32}\}$, then by the Gronwall inequality, we obtain
\begin{equation}\label{3.100-6}
\begin{aligned}
 \left\|\Delta u_{2}\right\|^2& \leq e^{-\int_{t_0-\tau}^t c_{30} d s} u_2\left(t_0-\tau\right)+c_{33}\int_{t_0-\tau}^t e^{\int_t^s c_{30} d y} \|h(x, s)\|^2 d s \\
& =c_{33}\int_{t_0-\tau}^t e^{-c_{30}(t-s)}\|h(x, s)\|^2 d s,\quad\forall \,t \in\left[t_0-\tau, t_0\right] .
\end{aligned}
\end{equation}

Thanks to \eqref{3.60-6} and \eqref{3.100-6}, we deduce there exists a constant $W_1>0$ depending on $c_{30},c_{31},c_{32},c_{33},t_0,\tau$ and $\|\widetilde{D}_{\delta}(t_0-\tau)\|_{X_t}^2$ such that
\begin{equation}\label{3.101-6}
\left\|\Delta u_2\right\|=\left\|\Delta U_{a_{2, \delta}}\left(t, t_0-\tau\right) u_{t_0-\tau}^0\right\|^2 \leq W_1<+\infty
\end{equation}
for any $t\in[t_0-\tau,t_0]$ and $\left(u_{t_0-\tau}^0 ; u_{t_0-\tau}^{1}\right) \in{\widetilde{D}_{\delta}}\left(t_0-\tau\right)$.

From \eqref{3.42-6} and $U_{b, \delta}\left(t, t_0-\tau\right) u_{t_0-\tau}^1) =  \partial_t u(t)$, it follows that there exists a constant $W_2>0$ depending on $\tau,T$ and $Z(\cdot)$ such that
\begin{equation}\label{3.102-6}
\left\|\nabla U_{b, \delta}\left(t, t_0-\tau\right) u_{t_0-\tau}^1\right\|^2 \leq W_2<+\infty
\end{equation}
for any $\tau>0$ and any $\left(u_{t_0-\tau}^0 ; u_{t_0-\tau}^{1}\right) \in{\widetilde{D}_{\delta}}\left(t_0-\tau\right)$.

Hence, from Lemma \ref{lem2.8-6}, we deduce the process $U_{\delta}(\cdot,\cdot)$ generated by problem \eqref{1.1-6} is pullback $\mathcal{D_\delta}$-asymptotically compact in $X_t$.

Then by Lemma \ref{lem2.5-6}, we derive there exists a pullback attractor $\widetilde{\mathcal A}_\delta=\{\widetilde{\mathcal A}(t)\}_{t \in \mathbb{R}}$ in $X_t$, which satisfies \eqref{3.85-6}. $\hfill$$\Box$

\begin{Remark}\label{remar3.8-6}
From Lemmas $\ref{lem3.3-6}$ and $\ref{lem3.5-6}$ and Theorem $\ref{th3.7-6}$, then it follows that the process $U_{\delta}(\cdot,\cdot)$ generated by problem \eqref{1.1-6} is pullback $\mathcal{D_\delta}$-asymptotically compact in $X_t$ and $\widetilde{ A}_\delta(t)=A_\delta(t):=\Lambda\left(D_\delta, t\right)$ for any $t\in\mathbb R$ and any $\delta\geq 0$.
\end{Remark}

\section{Upper semicontinuity of pullback attractors}
\ \ \ \ In this section, by using Lemmas \ref{lem2.9-6} and \ref{lem2.10-6}, we shall investigate the relationship between $A_\delta(t)$ and $A_0(t)$, which are pullback attractors of processes generated by problem \eqref{1.1-6} for $\delta>0$ and $\delta=0$ in the time-dependent space $X_t$, respectively.

The following theorem is the main result of this section.

\begin{Theorem}
Under the assumptions of Theorem $\ref{th3.7-6}$, then the pullback attractor $\mathcal{A}_{\delta}=\{A_\delta(t)\}_{t\in\mathbb R}$ given by Remark $\ref{remar3.8-6}$ satisfies
\begin{equation}\label{4.1-6}
\lim _{\delta \rightarrow 0^+} \sup_{t \in\left[\alpha, \beta\right]} \operatorname{dist}_{X_t}\left(A_\delta(t), A_0(t)\right)=0
\end{equation}
for any $[\alpha, \beta]\subset \mathbb R$ and $\delta\geq0$.
\label{th4.1-6}
\end{Theorem}

In order to prove Theorem \ref{th4.1-6}, it is necessary to establish the following three lemmas, then from Lemmas \ref{lem2.9-6} and \ref{lem2.10-6}, we can easily derive Theorem \ref{th4.1-6}.

\begin{Lemma} \label{lem4.2-6}
Under the assumptions of Theorem $\ref{th4.1-6}$, assume the family $\widetilde {\mathcal{D}}_{\delta}=\{\widetilde {\mathcal{D}}_{\delta}(t)\}_{t \in \mathbb{R}}$ satisfies \eqref{3.66-6}, then there exists a $\delta\in[0,1]$ such that
\begin{equation}\label{4.2-6}
{\overline {\bigcup\limits_{\delta  \in [0,1]} {{D_\delta }(t)} } ^{X_t}} \subset {D_{{0}}}(t)
\end{equation}
for any $t\in \mathbb R$.
\end{Lemma}
$\mathbf{Proof.}$ By Lemmas \ref{lem3.3-6} and \ref{lem3.5-6}, we conclude there exists a $\delta \in[0,1]$ such that \eqref{4.2-6} holds. $\hfill$$\Box$

\begin{Lemma} \label{lem4.3-6}
Under the assumptions of Theorem $\ref{th4.1-6}$, assume $\tau>0$, $\epsilon_i\rightarrow0^+$, $x_i\rightarrow x_0$ in the time-dependent space $X_t$ and $t_i>\alpha$ with $t_i\rightarrow t_0$ for any $i\in\mathbb N$, then
\begin{equation}\label{4.3-6}
U_{\epsilon_i}\left(t_i, \alpha-\tau\right) x_i \rightarrow U_0\left(t_0, \alpha-\tau\right) x_0
\end{equation}
in the time-dependent space $X_t$.
\end{Lemma}
$\mathbf{Proof.}$ Suppose $u_i(t)=U_{\delta_i}(t, \alpha-\tau) x_i$ is a weak solution to problem \eqref{1.1-6} with respect to $\delta=\delta_i$ and initial value $\left(u_{\alpha-\tau}^0 ; u_{\alpha-\tau}^1\right)=(u_{\alpha-\tau,i}^0 ; u_{\alpha-\tau, i}^{1})=x_i$.

Similarly, let $u_0=U_{0}(t, \alpha-\tau) x_0$ be a weak solution to problem \eqref{1.1-6} with respect to $\delta=0$ and initial value $\left(u_{\alpha-\tau}^0 ; u_{\alpha-\tau}^1\right)=\left(u_{\alpha-\tau, 0}^0 ; u_{\alpha-\tau, 0}^{1}\right)=x_0$.

Assume $t_i\in(\alpha,T)$ and
\begin{equation}\label{4.4-6}
\left\|\nabla u_{\alpha-\tau, i}^0\right\|^2+\varepsilon(\alpha-\tau)\left\|u_{\alpha-\tau, i}^1\right\|^2 \leq V_1, \quad\forall\,i\in\mathbb N,
\end{equation}
where $V_1>0$ is a constant.

From \eqref{3.86-6}, we deduce there exists a constant $V_2=V_2(\alpha-\tau,T,V_1)>0$ such that
\begin{equation}\label{4.5-6}
\left\|\nabla u_i(t)\right\|^2+\int_{\alpha-\tau}^T\left\|\nabla \partial_t u_i(s)\right\|^2 d s \leq V_2,\quad \forall\, t \in[\alpha-\tau, T], \,i \in \mathbb N.
\end{equation}

Let $z=u_a-u_b$ and $\delta=\delta_a-\delta_b$, then $z$ satisfies
\begin{equation}\label{4.6-6}
\begin{aligned}
& \varepsilon (t) \partial_{t t} z-\Delta z-\delta_a\left\|\nabla u_a\right\|^2 \Delta u_a+\delta_b\left\|\nabla u_b\right\|^2 \Delta u_b-\Delta \partial_t z+\lambda z \\
& =g\left(u_a\right)-g\left(u_b\right).
\end{aligned}
\end{equation}

Then inserting
\begin{align}\label{4.7-6}
& -\frac{1}{2}\left(\delta_a\left\|\nabla u_a\right\|^2+\delta_b\left\|\nabla u_b\right\|^2\right) \Delta z-\frac{1}{2} \delta_a\left(\left\|\nabla u_a\right\|^2-\left\|\nabla u_b\right\|^2\right) \Delta\left(u_a+u_b\right) \non\\
& -\frac{1}{2}\left(\delta_a-\delta_b\right)\left\|\nabla u_b\right\|^2 \Delta(u_a+u_b)\non \\
& =-\delta_a\left\|\nabla u_a\right\|^2 \Delta u_a+\delta_b\left\|\nabla u_b\right\|^2 \Delta u_b
\end{align}
into \eqref{4.6-6}, we obtain
\begin{align}\label{4.8-6}
& \varepsilon(t) \partial_{t t} z-\Delta z-\frac{1}{2}\left(\delta_a\left\|\nabla u_a\right\|^2+\delta_b\left\|\nabla u_b\right\|^2\right) \Delta z-\frac{1}{2} \delta_a\left(\left\|\nabla u_a\right\|^2-\left\|\nabla u_b\right\|^2\right) \Delta\left(u_a+u_b\right) \non\\
& -\frac{1}{2}\left(\delta_a-\delta_b\right)\left\|\nabla u_b\right\|^2 \Delta\left(u_a+u_b\right)-\Delta \partial_t z+\lambda z \non\\
& =g\left(u_a\right)-g\left(u_b\right) .
\end{align}

Additionally, the initial data of $z$ is
\begin{equation}\label{4.9-6}
\left(z_{\alpha-\tau}^0 ; z_{\alpha-\tau}^{1}\right)=\left(u_{\alpha-\tau, a}^0-u_{\alpha-\tau, b}^0 ; u_{\alpha-\tau, a}^{\prime}-u_{\alpha-\tau, b}^{\prime}\right).
\end{equation}

Choosing $\partial_t z+\xi z$ as a test function to \eqref{4.8-6}, we conclude
\begin{align}\label{4.10-6}
& \frac{1}{2} \frac{d}{d t}\left(\varepsilon(t)\|\partial_tz\|^2+2 \xi \varepsilon(t)(\partial_t z, z)+\|\nabla z\|^2+\xi\|\nabla z\|^2+\lambda\|z\|^2\right)-\frac{1}{2} \varepsilon^{\prime}(t)\|\partial_t z\|^2 \non\\
& -\xi \varepsilon^{\prime}(t)(\partial_tz, z)-\xi \varepsilon(t)\|\partial_t z\|^2+\xi\|\nabla z\|^2 +\frac{\xi}{2}\left(\delta_a\left\|\nabla u_a\right\|^2+\delta_b\left\|\nabla u_b\right\|^2\right)\|\nabla z\|^2\non\\
& +\frac{\xi \delta_a}{2}\left(\left\|\nabla u_a\right\|^2-\left\|\nabla u_b\right\|^2\right)^2 +\frac{1}{2}\left(\delta_a\left\|\nabla u_a\right\|^2+\delta_b\left\|\nabla u_b\right\|^2\right)(-\Delta z, \partial_t z)\non\\
& +\frac{\delta_a}{2}\left(\left\|\nabla u_a\right\|^2-\left\|\nabla u_b\right\|^2\right)\left(-\Delta\left(u_a+u_b\right), \partial_tz\right)  +\|\nabla \partial_tz\|^2+\lambda \xi\|z\|^2\non\\
& +\frac{\delta_a-\delta_b}{2}\left\|\nabla u_b\right\|^2\left(-\Delta\left(u_a+u_b\right), \partial_tz\right)+\frac{\xi\left(\delta_a-\delta_b\right)}{2}\left\|\nabla u_b\right\|^2\left(\left\|\nabla u_a\right\|^2-\left\|\nabla u_b\right\|^2\right) \non\\
&=\left(g\left(u_a\right)-g\left(u_b\right), \partial_t z+\xi z\right),
\end{align}
where $\xi>0$ is a constant.

Furthermore, let
\begin{align}\label{4.11-6}
\qquad\qquad\quad\,\, E_0{(t)}=&-\frac{1}{2} \varepsilon^{\prime}(t)\|\partial_tz\|^2-\xi \varepsilon^{\prime}(t)(\partial_t z, z)-\xi \varepsilon(t)\|\partial_t z\|^2+\xi\|\nabla z\|^2 \non\\
&+\frac{\xi}{2}\left(\delta_a\left\|\nabla u_a\right\|^2+\delta_b\left\|\nabla u_b\right\|^2\right)\|\nabla z\|^2+\|\nabla \partial_t z\|^2+\lambda \xi\|z\|^2 \non\\
&+\frac{\xi \delta_a}{2}\left(\left\|\nabla u_a\right\|^2-\left\|\nabla u_b\right\|^2\right)^2,
\end{align}
\begin{equation}\label{4.12-6}
\begin{aligned}
\quad\, E_1(t)= & \,\frac{1}{2}\left(\delta_a\left\|\nabla u_a\right\|^2+\delta_b\left\|\nabla u_b\right\|^2\right)(-\Delta z, \partial_t z)\\
&+\frac{\delta_a}{2}\left(\left\|\nabla u_a\right\|^2-\left\|\nabla u_b\right\|^2\right)\left(-\Delta\left(u_a+u_b\right), \partial_t z\right),
\end{aligned}
\end{equation}
\begin{equation}\label{4.13-6}
\begin{aligned}
E_2(t)= \,&\frac{\delta_a-\delta_b}{2}\left\|\nabla u_b\right\|^2\left(-\Delta\left(u_a+u_b\right), \partial_tz\right) \qquad \qquad\\
&+\frac{\xi\left(\delta_a-\delta_b\right)}{2}\left\|\nabla u_b\right\|^2\left(\left\|\nabla u_a\right\|^2-\left\|\nabla u_b\right\|^2\right)
\end{aligned}
\end{equation}
and
\begin{equation}\label{4.14-6}
E_3(t)=\left(g\left(u_a\right)-g\left(u_b\right), \partial_t z+\xi z\right).\qquad \qquad \qquad\quad
\end{equation}

Then inserting $\eqref{4.11-6}-\eqref{4.14-6}$ into \eqref{4.10-6}, we derive
\begin{align}\label{4.15-6}
\frac{1}{2} \frac{d}{d t}& \left(\varepsilon(t)\left\|\partial_t z\right\|^2+2 \xi \varepsilon(t)\left(\partial_t z, z\right)+\|\nabla z\|^2+\xi\|\nabla z\|^2+\lambda\|z\|^2\right) \non\\
&+E_0(t)+E_1(t)+E_2(t)=E_3(t).
\end{align}

By \eqref{4.5-6}, the Cauchy and Young inequalities, we deduce
\begin{align}\label{4.16-6}
\left|E_1(t)\right| & \leq \frac{1}{2}\left(\delta_a\left\|\nabla u_a\right\|^2+\delta_b\left\|\nabla u_b\right\|^2\right)\left(\frac{1}{2}\|\nabla z\|^2+\frac{1}{2}\left\|\nabla \partial_t z\right\|^2\right) \non\\
& +\left|\frac{\delta_a}{2}\left(\left\|\nabla u_a\right\|^2-\left\|\nabla u_b\right\|^2\right)\left[\left(\nabla u_a, \nabla \partial_t z\right)+\left(\nabla u_b, \nabla \partial_t z\right)\right]\right|\non \\
& \leq \frac{1}{4}\left(\delta_a+\delta_b\right) V_2^2+\frac{1}{4}\left(\delta_a+\delta_b\right) V_2\left\|\nabla \partial_t z\right\|^2 \non\\
& +\frac{\delta_a}{2}\left|\left\|\nabla u_a\right\|^2\left[\left(\nabla u_a, \nabla \partial_t z\right)+\left(\nabla u_b, \nabla \partial_t z\right)\right]\right|\non\\
& \leq \frac{3}{4} \delta_a V_2^2+\frac{1}{4} \delta_b V_2^2+\frac{3}{4} \delta_a V_2\|\nabla \partial_t z\|^2+\frac{1}{4} \delta_b V_2\|\nabla \partial_t z\|^2 \non\\
& \left.\leq2(\delta_a+\delta_b\right) V_2^2+2\left(\delta_a+\delta_b\right) V_2\left\|\nabla \partial_t z\right\|^2
\end{align}
and
\begin{align}\label{4.17-6}
\left|E_2(t)\right| & \leq \frac{\left|\delta_a-\delta_b\right|}{32 \xi}\left\|\nabla u_b\right\|^2\left\|\nabla u_a\right\|^2+\frac{\left|\delta_a-\delta_b\right|}{4}\left\|\nabla u_b\right\|^2\left\|\nabla \partial_t z\right\|^2\non \\
& +\frac{\xi\left|\delta_a-\delta_b\right|}{2}\left\|\nabla u_b\right\|^2\left\|\nabla u_a\right\|^2\non \\
& \leq\left|\delta_a-\delta_b\right|(1+\delta) V_2^2+\left|\delta_a-\delta_b\right|V_2\left\|\nabla \partial_t z\right\|^2,
\end{align}
where $t\in[\alpha-\tau,T]$ and $a,b\in\mathbb N$.

Using \eqref{1.4-6}, \eqref{1.5-6}, the Young inequality and the embedding $H_0^1 (\Omega)\hookrightarrow L^{\frac{2 n}{n-2}}(\Omega)$, we derive there exists a constant $C=C(k,\xi,V_2)>0$ such that
\begin{equation}\label{4.18-6}
E_3(t) \leq C\left(\left\|\nabla z\right\|^2+\|\partial_t z\|^2\right)+\frac{1}{2}\|\nabla \partial_t z\|^2.
\end{equation}

Let
\begin{equation}\label{4.19-6}
\widetilde E(t)=\varepsilon(t)\|\partial_tz\|^2+2 \xi \varepsilon(t)(\partial_t z, z)+\|\nabla z\|^2+\xi\|\nabla z\|^2+\lambda\|z\|^2,
\end{equation}
then from $\eqref{4.15-6}-\eqref{4.19-6}$, we obtain there exists a constant $\widetilde C=\widetilde C(\delta,k,L,V_2)>0$ such that
\begin{equation}\label{4.20-6}
\frac{d}{d t} \widetilde{E}(t) \leq \widetilde{C} \widetilde{E}(t)+\widetilde{C}\left(\delta_a+\delta_b\right)\left(\|\nabla \partial_t z(t)\|^2+1\right).
\end{equation}

Then by the Gronwall inequality, we arrive at
\begin{equation}\label{4.21-6}
\begin{aligned}
\widetilde{E}(t) & \leq  e^{ \widetilde C(T-\alpha+\tau)} \widetilde{E}(\alpha-\tau)+ \widetilde C\left(\delta_a+\delta_b\right) \int_{\alpha-\tau}^T e^{ \widetilde C(T-s)}\left(\|\nabla \partial_tz(s)\|^2+1\right) d s.
\end{aligned}
\end{equation}

Thanks to \eqref{4.5-6} and \eqref{4.19-6}, then there exists a constant $C>0$ such that
\begin{equation}\label{4.22-6}
\widetilde C\left(\delta_a+\delta_b\right) \int_{\alpha-\tau}^T e^{ \widetilde C(T-s)}\left(\|\nabla \partial_tz(s)\|^2+1\right) d s \leq C.
\end{equation}

Moreover, noticing $u_i(t)=U_{\delta_i}(t, \alpha-\tau) x_i$ and $z=u_a-u_b$, we derive
\begin{equation}\label{4.23-6}
u_a(t)=U_{\delta_a}(t, \alpha-\tau) x_a
\end{equation}
and
\begin{equation}\label{4.24-6}
u_b(t)=U_{\delta_b}(t, \alpha-\tau) x_b.
\end{equation}

From \eqref{4.19-6}, \eqref{4.21-6}, \eqref{4.23-6} and \eqref{4.24-6}, we deduce
\begin{equation}\label{4.25-6}
\begin{aligned}
\widetilde E(t) &= \left\|U_{\delta_a}(t, \alpha-\tau) x_a-U_{\delta_b}(t, \alpha-\tau) x_b\right\|_{X_t}^2+\xi(\|\nabla z(t)\|^2+2 \varepsilon(t)(\partial_t z(t), z(t))) \\
& +\lambda\|z(t)\|^2.
\end{aligned}
\end{equation}

From \eqref{4.21-6} and \eqref{4.25-6}, we obtain
\begin{equation}\label{4.26-6}
\begin{aligned}
&\left\|U_{\delta_a}(t, \alpha-\tau) x_a-U_{\delta_b}(t, \alpha-\tau) x_b\right\|_{X_t}^2+\xi(\|\nabla z(t)\|^2+2 \varepsilon(t)(\partial_t z(t), z(t))+\lambda\|z(t)\|^2 \\
& \leq  e^{ \widetilde C(T-\alpha+\tau)} \widetilde{E}(\alpha-\tau)+ \widetilde C\left(\delta_a+\delta_b\right) \int_{\alpha-\tau}^T e^{ \widetilde C(T-s)}\left(\|\nabla \partial_tz(s)\|^2+1\right) d s.
\end{aligned}
\end{equation}

Suppose $\widetilde {C}_1\geq \max\{2\xi L V_2,\widetilde C\}$, then by \eqref{4.5-6} and \eqref{4.26-6}, we conclude
\begin{align}\label{4.27-6}
& \left\|U_{\delta_a}(t, \alpha-\tau) x_a-U_{\delta_b}(t, \alpha-\tau) x_b\right\|_{X_t}^2 \non\\
& \leq e^{\widetilde{C}_1(T-\alpha+\tau)}\left[\varepsilon(\alpha-\tau)\left\|u_{\alpha-\tau, a}^{1}-u_{\alpha-\tau, b}^1\right\|^2\right.\non\\
& \qquad\qquad\quad\quad+2 \xi \varepsilon(\alpha-\tau)\left(u_{\alpha-\tau, a}^1-u_{\alpha-\tau, b}^{1}, u_{\alpha-\tau, a}^0-u_{\alpha-\tau, b}^0\right) \non\\
& \qquad\qquad\quad\quad\left.+(1+\xi) \left\|\nabla\left(u_{\alpha-\tau, a}^0-u_{\alpha-\tau, b}^0\right)\right\|^2+\lambda\left\|u_{\alpha-\tau, a}^0-u_{\alpha-\tau, b}^0\right\|^2\right] \non\\
& +\widetilde{C}_1\left(\delta_a+\delta_b\right) \int_{\alpha-\tau}^T e^{\widetilde{C}_1(T-s)}\left(\|\nabla \partial_t z(s)\|^2+1\right) d s.
\end{align}

Using \eqref{4.22-6} and \eqref{4.27-6}, we derive there exists a constant $V_3=V_3(\widetilde{C}_1, V_2,\xi,\lambda)>0$ such that
\begin{equation}\label{4.28-6}
\begin{aligned}
& \left\|U_{\delta_a}(t, \alpha-\tau) x_a-U_{\delta_b}(t, \alpha-\tau) x_b\right\|_{X_t}^2 \\
& \leq V_3\left(\left\|\nabla\left(u_{\alpha-\tau, a}^0-u_{\alpha-\tau, b}^0\right)\right\|^2+\left\|u_{\alpha-\tau, a}^1-u_{\alpha-\tau, b}^1\right\|^2+\delta_a+\delta_b\right)
\end{aligned}
\end{equation}
for any $t\in[\alpha-\tau,T]$ and $a,b\in\mathbb N$.

Furthermore, it follows from \eqref{4.28-6} that there exists a $\mathbb N_{\delta}\in\mathbb N$ such that
\begin{equation}\label{4.29-6}
\left\|U_{\delta_a}(t, \alpha-\tau) x_a-U_{\delta_b}(t, \alpha-\tau) x_b\right\|_{X_t}^2 \leq \eta,
\end{equation}
where $\eta>0$ is a constant, $t\in[\alpha,T]$ and $a,b\in\mathbb N$.

Therefore, we obtain $\left\{U_{\delta_i}(t, \alpha-\tau) x_i\right\}_{i \in \mathbb N}$ is a Cauchy sequence in $X_t$ for any $t \in[\alpha, T]$.

Without loss of generality, let $\left\{U_{\delta_i}(\tilde t, \alpha-\tau) x_i\right\}_{i \in \mathbb N}$ be a convergent sequence in $X_t$, that is, for some $\tilde y\in X_t$,
\begin{equation}\label{4.30-6}
U_{\delta_i}(\tilde{t}, \alpha-\tau) x_i \rightarrow \tilde{y} \quad\text { as } i \rightarrow+\infty.
\end{equation}

Then it follows that
\begin{align}\label{4.31-6}
&\left\|U_{\delta_i}\left(t_i, \alpha-\tau\right) x_i-U_{\delta_i}(\tilde{t}, \alpha-\tau) x_i\right\|_{X_t}\non \\
& = \| U_{\delta_i}\left(t_i, \alpha-\tau\right) x_i-U_{\delta_i}(\tilde{t}, \alpha-\tau) x_i-U_{\delta_b}(t_i,\alpha-\tau)x_b+ U_{\delta_b}(t_i,\alpha-\tau)x_b\non\\
 & -U_{\delta_b}(\tilde{t}, \alpha-\tau) x_b+U_{\delta_b}(\tilde{t}, \alpha-\tau) x_b \|_{X_t} \non\\
 &\leq \left\|U_{\delta_i}\left(t_i, \alpha-\tau\right) x_i-U_{\delta_b}\left(t_i, \alpha-\tau\right) x_b\right\|_{X_t} +\left\|U_{\delta_b}(t_i, \alpha-\tau) x_b-U_{\delta_b}(\tilde{t}, \alpha-\tau) x_b\right\|_{X_t}\non\\
& +\left\|U_{\delta_b}(\tilde{t}, \alpha-\tau) x_b-U_{\delta_i}(\tilde{t}, \alpha-\tau) x_i\right\|_{X_t}.
\end{align}

Next, we will take turns to analyze each item in \eqref{4.31-6}.

Noticing $t_i\rightarrow t_0$ as $i\rightarrow+\infty$, then we can find suitable $\tilde b,n_1\in\mathbb N$ such that
\begin{align}\label{4.32-6}
&\,\quad\left\|U_{\delta_i}\left(t_i, \alpha-\tau\right) x_i-U_{\delta_i}(\tilde{t}, \alpha-\tau) x_i\right\|_{X_t} \non\\
&\leq \left\|U_{\delta_i}\left(t_i, \alpha-\tau\right) x_i-U_{\delta_{\tilde{b}}}\left(t_i, \alpha-\tau\right) x_{\tilde{b}}\right\|_{X_t}+\|U_{\tilde{\delta}_{\tilde{b}}}\left(t_i, \alpha-\tau\right) x_{\tilde{b}}-U_{\delta_i}(\tilde{t}, \alpha-\tau) x_i \|_{X_t} \non\\
&\leq  \left\|U_{\delta_i}\left(t_i, \alpha-\tau\right) x_i-U_{\delta_{\tilde{b}}}\left(t_i, \alpha-\tau\right) x_{\tilde{b}}\right\|_{X_t}+\left\|U_{\delta_{\tilde{b}}}(\tilde{t}, \alpha-\tau) x_{\tilde{b}}-U_{\delta_i}(\tilde{t}, \alpha-\tau) x_i\right\|_{X_t}\non\\
&\leq \frac{\eta}{2}.
\end{align}

By Theorem \ref{th3.1-6} and $u_i(t)=U_{\delta_i}(t, \alpha-\tau) x_i$, we deduce $U_{\delta_{\widetilde{k}}}(\cdot, \alpha-\tau) \in C\left([\alpha-\tau, T] ; X_t\right)$, which combines with $t_i\rightarrow\tilde t$ as $i\rightarrow +\infty$ to yield that there exists a $n_2\geq n_1$ such that
\begin{equation}\label{4.33-6}
\left\|U_{\delta_{\tilde{b}}}\left(t_i, \alpha-\tau\right) x_{\tilde{b}}-U_{\delta_{\tilde{b}}}(\tilde{t}, \alpha-\tau) x_{\tilde{b}}\right\|_{X_t} \leq \frac{\eta}{2}
\end{equation}
when $n\geq n_2$.

From $\eqref{4.31-6}-\eqref{4.33-6}$, we obtain there exists a $n_{\eta}\in \mathbb N$ large enough such that
\begin{equation}\label{4.34-6}
\left\|U_{\delta_i}\left(t_i, \alpha-\tau\right) x_i-U_{\delta_i}(\tilde{t}, \alpha-\tau) x_i\right\|_{X_t} \leq \eta, \quad \forall\, n \geq n_\eta.
\end{equation}

Thus by \eqref{4.30-6} and \eqref{4.34-6}, we arrive at
\begin{equation}\label{4.35-6}
U_{\delta_i}\left(t_i, \alpha-\tau\right) x_i \rightarrow \tilde{y}  \quad\text { as } i \rightarrow+\infty
\end{equation}
in $X_t$.

In addition, if let the constant $V_2$ in \eqref{4.5-6} be large enough, then \eqref{4.34-6} also holds for $u_0$.

Suppose $v=u_i-u_0$, then $v$ satisfies
\begin{equation}\label{4.36-6}
\varepsilon {(t)} \partial_{t t} v-\Delta v-\delta_i\left\|\nabla u_i\right\|^2 \Delta u_i-\Delta \partial_t v+\lambda v=g\left(u_i\right)-g\left(u_0\right)
\end{equation}
with initial data
\begin{equation}\label{4.37-6}
\left(v_{\alpha-\tau}^0 ; v_{\alpha-\tau}^1\right)=\left(u_{\alpha-\tau, i}^0-u_{\alpha-\tau, 0}^0 ;u_{\alpha-\tau, i}^{1}-u_{\alpha-\tau, 0}^1\right).
\end{equation}

From Lemma \ref{lem3.2-6}, we obtain
\begin{equation}\label{4.38-6}
\left\|\nabla\left(U_{\delta_i}(\tilde{t}, \alpha-\tau) x_i-U_0(\tilde{t}, \alpha-\tau) x_0\right)\right\|^2+\varepsilon(\tilde t)\left\|\partial_t\left(U_{\delta_i}(\tilde{t}, \alpha-\tau) x_i-U_0(\tilde{t}, \alpha-\tau) x_0\right)\right\|^2 \rightarrow 0.
\end{equation}

Then by \eqref{4.30-6} and \eqref{4.35-6}, we conclude
\begin{equation}\label{4.39-6}
U_{\delta_i}\left(t_i, \alpha-\tau\right) x_i \rightarrow U_0\left(t_0, \alpha-\tau\right) x_0
\end{equation}
in $X_t$.
$\hfill$$\Box$

Inspired by Lemmas \ref{lem2.9-6} and \ref{lem2.10-6}, if we can verify the three properties (i)$-$(iii) in Lemma \ref{lem2.10-6}, then the following lemma holds immediately. To this end, we will perform some calculations similar to those in the proof of Theorem \ref{th3.7-6}.

\begin{Lemma} \label{lem4.4-6}
Under the assumptions of Theorem $\ref{th4.1-6}$, then for any $\delta_i \in[0,1],t_i \in[\alpha, \beta],\tau_i \rightarrow+\infty$ and $x_i \in D_{\delta_ i}\left(t_i-\tau_i\right)$ with $i\in\mathbb N$, the sequence $\left\{U_{\delta_i}\left(t_i, t_i-\tau_i\right) x_i\right\}_{i \in \mathbb N}$ is relatively compact in $X_t$.
\end{Lemma}
$\mathbf{Proof.}$ Thanks to \eqref{3.94-6}, we deduce
\begin{equation}\label{4.40-6}
\left\|\nabla v_{a_1, \delta}\left(t, t_0-\tau\right) u^0_{t_0-\tau}\right\|^2 \leq e^{-2 \tau}\left\|\nabla u^0_{t_0-\tau}\right\|^2,
\end{equation}
which leads to
\begin{equation}\label{4.41-6}
\sup _{\delta \in[0, 1]}\left\|\nabla U_{a_1, \delta}(t, \alpha-\tau) u_{\alpha-\tau}^0\right\|^2 \leq e^{-2 \tau}\left\|\nabla u_{\alpha-\tau}^0\right\|^2.
\end{equation}

From \eqref{3.101-6}, we obtain there exists a constant $\widetilde {W}_1>0$ such that
\begin{equation}\label{4.42-6}
\sup _{\delta \in[0, 1]}\left\|\Delta U_{a_2, \delta}(t, \alpha-\tau) u_{\alpha-\tau}^0\right\|^2 \leq  \widetilde {W}_1.
\end{equation}

Furthermore, by \eqref{3.102-6}, we deduce there exists a constant $\widetilde {W}_2>0$ such that
\begin{equation}\label{4.43-6}
\sup _{\delta[0,1]}\left\|\nabla U_{b, \delta}(t, \alpha-\tau) u_{\alpha-\tau}^{1}\right\|^2 \leq \widetilde{W}_2.
\end{equation}

The above inequalities $\eqref{4.40-6}-\eqref{4.42-6}$ hold for any $t\in[\alpha,\beta]$, $\tau>0$ and $\left(u_{\alpha-\tau}^0 ; u_{\alpha-\tau}^{1}\right) \in D_\delta(t-\tau)$, $\widetilde{W}_1$ and $\widetilde{W}_2$ depend on $\alpha-\tau$ and $\left\|D_\delta(\alpha-\tau)\right\|_{X_t}$.

As a result, by Lemmas \ref{lem2.9-6} and \ref{lem2.10-6}, then it follows that $\left\{U_{\delta_i}\left(t_i, t_i-\tau_i\right) x_i\right\}_{i \in \mathbb N}$ is relatively compact in $X_t$.
$\hfill$$\Box$

$\mathbf{Acknowledgment}$

We would like to express our sincere gratitude to the editors and anonymous reviewers for their invaluable feedback and insightful comments on the initial draft of this paper. Their thorough review and constructive suggestions have played a crucial role in improving the quality of this manuscript. We greatly appreciate their dedication and support throughout the review process.
\newpage
$\mathbf{Funding}$

Yang was supported by the China Scholarship Council with number 202206630048. Qin was supported by the National Natural Science Foundation of China with contract numbers 12171082 and 12226403, the fundamental research funds for the central universities with contract numbers $2232022G$-$13$, $2232023G$-$13$ and a grant from science and technology commission of Shanghai municipality. Wang was supported by Science and Technology Commission of Shanghai Municipality with contract numbers 23ZR1402100.

$\mathbf{Conflict\,\,of\,\,interest\,\,statement}$

The authors have no conflict of interest.


\newpage
\begin{thebibliography}{lllp}
\setlength{\itemsep}{- 2mm}
\bibitem{bv1992} A. V. Babin and M. I. Vishik, Attractors of evolutionary equations, North-Holland, Amsterdam, 1992.




\bibitem{clr.6}T. Caraballo, G. Lukaszewicz and J. Real, Pullback attractors for asymptotically compact non¨Cautonomous dynamical systems, Nonlinear Anal., 64(3)(2006), 484-498.

\bibitem{c.5} I. Chueshov, Long-time dynamics of Kirchhoff wave models with strong nonlinear damping, J. Diff. Equ., 252(1)(2012), 1229-1262.

\bibitem{ConPT.6} M. Conti, V. Pata and R. Temam, Attractors for processes on time-dependent spaces, Applications to wave equations, J. Diff. Equ., 255(6)(2013), 1254-1277.




\bibitem{fz.6} X. Fan and Z. Zhou, Kernel sections for non-autonomous strongly damped wave equations of non-degenerate Kirchhoff-type, Appl. Mathe. Comput., 158(2)(2004), 253-266.




\bibitem{k.5} G. Kirchhoff, Vorlesungen $\rm \ddot{u}$ber Mechanik, Teubner, Sluttgart, 1883.




\bibitem{l.5} K. Li, A Gronwall-type lemma with parameter and its application to Kirchhoff type nonlinear wave equation, J. Math. Anal. Appl., 447(2)(2017), 683-704.

\bibitem{ly.5} Y. Li and Z. Yang, Robustness of attractors for non-autonomous Kirchhoff wave models with strong nonlinear damping, Appl. Math. Optim., 81(1)(2019), 245-272.

\bibitem{lyf.5} Y. Li, Z. Yang and N. Feng, Uniform attractors and their continuity for the non-autonomous Kirchhoff wave models, Disc. Contin. Dynam. Syst., 26(12)(2021), 6267-6284.

\bibitem{ljs.5} G. Liu and M. A. Jorge Silva, Attractors and their properties for a class of Kirchhoff models with integro-differential damping, Appl. Analy., 101(9)(2020), 3284-3307.



\bibitem{MWX2021-Y0} H. Ma, J. Wang and J. Xie, Pullback attractors for nonautonomous degenerate Kirchhoff equations  with strong damping, Adv. Math. Phy., (2021), 7575078, DOI: 10.1155/2021/7575078.

\bibitem{mz.5} H. Ma and C. Zhong, Attractors for the Kirchhoff equations with strong nonlinear damping, Appl. Math. Lett., 7(1)(2017), 127-133.




\bibitem{psz} X. Peng, Y. Shang and X. Zheng, Pullback attractors of nonautonomous nonclassical diffusion equations with nonlocal diffusion, Zeitschrift f\"ur Angewandte Mathematik und Physik, 69(4)(2018), 1-14.




\bibitem{r} J. C. Robinson, Infinite-dimensional dynamical systems, Cambridge University Press, England, 2011.




\bibitem{wyh1.6} Y. Wang, Pullback attractors for nonautonomous wave equations with critical exponent, Nonlinear Anal., TMA, 68(2)(2008), 365-376.

\bibitem{wz.5} Y. Wang and C. Zhong, Upper semicontinuity of pullback attractors for nonautonomous Kirchhoff wave models, Disc. Conti. Dyn. Sys., 33(7)(2013), 3189-3209.





\bibitem{YL.6} Z. Yang and Y. Li, Criteria on the existence and stability of pullback exponential attractors and their application to non-autonomous kirchhoff wave models, Disc. Contin. Dynam. Syst., 38(5)(2018), 2629-2653.

\bibitem{LYUPP.6} Z. Yang and Y. Li, Upper semicontinuity of pullback attractors for non-autonomous Kirchhoff wave equations, Disc. Contin. Dynam. Syst., 24(9)(2019), 4899-4912.




\bibitem{zs2} K. Zhu and C. Sun, Pullback attractors for nonclassical diffusion equations with delays, J. Math. Phy., 56(9)(2015), 1-20.







\end{thebibliography}
\end{document}